\documentclass[12pt]{article}

\usepackage{amscd,amsmath, amssymb, fancyhdr, epsfig}
\usepackage[matrix,arrow]{xy}

\numberwithin{equation}{section}

\newcommand{\version}{version 3.1,\ \   Aug. 10, 2012}

\def\eqref#1{(\ref{#1})}
\newcommand{\goth}{\mathfrak}

\newcommand{\arrow}{{\:\longrightarrow\:}}

\newcommand{\C}{{\Bbb C}}
\newcommand{\SLnH}{SL(n,{\Bbb H})}
\newcommand{\R}{{\Bbb R}}

\renewcommand{\H}{{\Bbb H}}
\newcommand{\6}{\partial}
\def\1{\sqrt{-1}\:}
\newcommand{\restrict}[1]{{\left|_{{\phantom{|}\!\!}_{#1}}\right.}}
\newcommand{\cntrct}                
{\hspace{2pt}\raisebox{1pt}{\text{$\lrcorner$}}\hspace{2pt}}

\makeatletter
\def\x@arrow{\DOTSB\Relbar}
\def\xlongequalsignfill@{\arrowfill@\x@arrow\Relbar\x@arrow}
\newcommand{\xlongequal}[2][]{%
        \ext@arrow 0099\xlongequalsignfill@{#1}{#2}}
\def\xlongrightarrowfill@{\arrowfill@\relbar\relbar\longrightarrow}
\newcommand{\xlongrightarrow}[2][]{%
        \ext@arrow 0099\xlongrightarrowfill@{#1}{#2}}
\makeatother

\renewcommand{\bar}{\overline}
\renewcommand{\phi}{\varphi}
\renewcommand{\epsilon}{\varepsilon}
\renewcommand{\geq}{\geqslant}
\renewcommand{\leq}{\leqslant}

\newcommand{\comass}{\operatorname{\sf comass}}
\newcommand{\Av}{\operatorname{\sf Av}}
\newcommand{\End}{\operatorname{End}}

\newcommand{\Vol}{\operatorname{Vol}}

\newcommand{\Hol}{\operatorname{Hol}}

\newcommand{\codim}{\operatorname{codim}}

\renewcommand{\Re}{\operatorname{Re}}

\newcounter{Mycounter}[section]
\newcounter{lemma}[section]
\renewcommand{\thelemma}{{Lemma \thesection.\arabic{lemma}}}
\newcommand{\lemma}{%
    \setcounter{lemma}{\value{Mycounter}}
    \refstepcounter{lemma}
    \stepcounter{Mycounter}
    {\noindent \bf \thelemma:\ }}

\newcounter{claim}[section]
\renewcommand{\theclaim}{{Claim \thesection.\arabic{claim}}}
\newcommand{\claim}{%
    \setcounter{claim}{\value{Mycounter}}
    \refstepcounter{claim}
    \stepcounter{Mycounter}
    {\noindent \bf \theclaim:\ }}

\newcounter{sublemma}[section]
\newcounter{corollary}[section]
\renewcommand{\thecorollary}{{Corollary \thesection.\arabic{corollary}}}
\newcommand{\corollary}{%
    \setcounter{corollary}{\value{Mycounter}}
    \refstepcounter{corollary}
    \stepcounter{Mycounter}
    {\noindent \bf \thecorollary:\ }}

\newcounter{theorem}[section]
\renewcommand{\thetheorem}{{Theorem \thesection.\arabic{theorem}}}
\newcommand{\theorem}{%
    \setcounter{theorem}{\value{Mycounter}}
    \refstepcounter{theorem}
    \stepcounter{Mycounter}
    {\noindent \bf \thetheorem:\ }}

\newcounter{conjecture}[section]
\newcounter{proposition}[section]
\renewcommand{\theproposition}
      {{Proposition \thesection.\arabic{proposition}}}
\newcommand{\proposition}{%
    \setcounter{proposition}{\value{Mycounter}}
    \refstepcounter{proposition}
    \stepcounter{Mycounter}
    {\noindent \bf \theproposition:\ }}

\newcounter{definition}[section]
\renewcommand{\thedefinition}
      {{Definition~\thesection.\arabic{definition}}}
\newcommand{\definition}{%
    \setcounter{definition}{\value{Mycounter}}
    \refstepcounter{definition}
    \stepcounter{Mycounter}
    {\noindent \bf \thedefinition:\ }}

\newcounter{example}[section]
\newcounter{remark}[section]
\renewcommand{\theremark}{{Remark \thesection.\arabic{remark}}}
\newcommand{\remark}{%
    \setcounter{remark}{\value{Mycounter}}
    \refstepcounter{remark}
    \stepcounter{Mycounter}
    {\noindent \bf \theremark:\ }}

\newcounter{problem}[section]
\newcounter{question}[section]
\makeatletter

\@addtoreset{equation}{section} \@addtoreset{footnote}{section} \makeatother

\def\blacksquare{\hbox{\vrule width 5pt height 5pt depth 0pt}}
\def\endproof{\blacksquare}

\begin{document}
\begin{center}
{\LARGE\bf
Calibrations in hyperk\"ahler geometry\\[4mm]
}

Gueo Grantcharov, Misha Verbitsky\footnote{Misha Verbitsky is
partially supported by RFBR grant 10-01-93113-NCNIL-a,
RFBR grant 09-01-00242-a, Simons-IUM fellowship,
Science Foundation of the SU-HSE award No. 10-09-0015 and AG Laboratory HSE,
RF government grant, ag. 11.G34.31.0023.}

\end{center}

{\small \hspace{0.10\linewidth}
\begin{minipage}[t]{0.85\linewidth}
{\bf Abstract} \\
We describe a family of calibrations arising naturally on a hyperk\"ahler manifold $M$. These calibrations
calibrate the holomorphic Lagrangian, holomorphic isotropic and holomorphic coisotropic subvarieties. When $M$
is an HKT (hyperk\"ahler with torsion) manifold with holonomy $SL(n, {\Bbb H})$, we construct another family of
calibrations $\Phi_i$, which calibrates holomorphic Lagrangian and holomorphic coisotropic subvarieties. The
calibrations $\Phi_i$ are (generally speaking) not
parallel with respect to any torsion-free connection on $M$.


\end{minipage}
}

\tableofcontents


\section{Introduction}


The theory of calibrations was developed by R. Harvey and B. Lawson in \cite{_Harvey_Lawson:Calibrated_}, and proved
to be very useful in describing the geometric structures associated with special holonomies. Since then calibrations have become a central notion in many geometric developments in string physics and M-theory.
Up to dimension 8, the calibrations are thoroughly studied and pretty much understood (\cite{_DHM:R^8_}), but in the higher dimensions, the
classification problem seems to be immense. Even in more special situations, such as in hyperk\"ahler geometry,
the problem of classification of natural\footnote{For hyperk\"ahler geometry, ``natural'' would mean
``$Sp(n)$-invariant''.} calibrations is unsolved.

\hfill

On a K\"ahler manifold, the normalized power of the K\"ahler form $\frac{\omega^p}{p!}$ is a
calibration. A subvariety is complex
analytic if and only if it is calibrated. This is actually
very easy to see, because a subspace $V\subset TM$ is a
face of $\frac{\omega^p}{p!}$ if and only if $V$ is
complex linear (this follows {}from the so-called
``Wirtinger inequalities'', see
e.g. \cite{_Harvey_Lawson:Calibrated_}).

\hfill

In this paper we study a family of calibrations which appear naturally in quaternionic geometry, and describe
the corresponding calibrated subvarieties. These calibrations are in many ways analogous to the powers of the
K\"ahler form. We define several new calibrations, for hyperk\"ahler, hypercomplex and HKT-geometry.
{}From the calibration-theoretic point of view, the last of these is most interesting, because it is (generally
speaking) not preserved by {\em any} torsionless connection on $M$. Some of these forms were considered previously in 
\cite{_Verbitsky:skoda.tex_,_AV:Calabi_,_Verbitsky:balanced_}.

\hfill

In hyperk\"ahler geometry, the role of a K\"ahler form is
played by a 4-form $\Theta:= \omega_I^2 + \omega_J^2 +
\omega_K^2$. In Section 5.2 we show that the normalized powers $\Theta^p$ are calibrations. It is easy to
see that $V\subset TM$ is a face of $\Theta$ if and only
if $V$ is a quaternionic subspace
(\ref{_triana_calibra_Theorem_}).

The corresponding calibrated subvarieties are those which are complex analytic with respect to $I$, $J$ and $K$.
Such subvarieties are called {\bf trianalytic}. In \cite{_Verbitsky:trianalyt_,_Verbitsky:desing_}, the theory
of trianalytic subvarieties was developed to some extent. It was shown that the trianalytic subvarieties admit a
canonical desingularizaton, which is hyperk\"ahler. Also it was shown that any complex analytic subvariety of
$(M,I)$ is trianalytic, if the complex structure $I$ is generic in its twistor family.

\hfill

Any homogeneous polynomial $P(x,y, z)$ of degree $p$ gives
a closed $2p$-form $P(\omega_I, \omega_J, \omega_K)$
on $M$, and (when the holonomy of $M$ is maximal) all
parallel differential forms on $M$ are obtained this way.
When $P(x,y,z) = \frac{x^p}{p!}$, it is a K\"ahler calibration, when
$P(x,y,z)=c_p(x^2+y^2+z^2)^p$, where
$c_p=\sum_{k=0}^p \frac{(p!)^2}{(k!)^2}(2k)!4^{p-k}$,
it is the trianalytic calibration defined above( Theorem 5.3).
It would be interesting to classify all
calibrations obtained this way.

The calibrations $\Psi_k$ and $\Phi_{n+k}$ we study in this paper are also polynomials
on $\omega_I, \omega_J, \omega_K$. These calibrations are
called {\bf holomorphic Lagrangian},
{\bf holomorphic isotropic} and {\bf holomorphic
  coisotropic calibrations}.
The form $\Psi_k$ is obtained as a
$(k,k)$-component of $\Re(\omega_I-\1\omega_K)^k$,
normalized in appropriate way, where $\omega_I-\1\omega_K$
is a holomorphic symplectic form on $(M,J)$, and the
$(k,k)$-part is taken with respect to the complex structure $I$. In \cite{_Verbitsky:skoda.tex_,_AV:Calabi_}
it was proven that this form is closed and weakly positive.

\hfill

We show in Section 5.4 that a subvariety $Z\subset M$ is calibrated by
$\Psi_k$ if and only if $Z$ is holomorphic Lagrangian in
$(M, I)$ (for $k=\frac 1 2 \dim_\C M$) and isotropic (for
$k<\frac 1 2 \dim_\C M$) (\ref{_complex_Lagrangian_subspaces_Proposition_},
\ref{_iso_cali_Proposition_}). Note that holomorphic Lagrangian calibrations
have been found previously in \cite{BrH} in dimension eight. 

In \cite{JF} a different holomorphic Lagrangian calibration in any dimension was constructed as part of an investigation relating the faces of some calibrations to intersecting supersymmetric branes in M-theory. In String Theory the holomorphic
Lagrangian submanifolds were related to 3-dimensional
topological field theory with target hyperk\"ahler
manifold \cite{KRS}. In Section 5.6 we provide some examples of holomorphic Lagrangian subvarieties of hypercomplex manifolds which are not hyperk\"ahler.

The proof of this result relies on a particular partial
order defined on the set of precalibrations. We say that
$\eta \preceq \eta_1$ if all faces of $\eta$ are also
faces of $\eta_1$. For instance, the calibrations
$c_p\Theta^p$, $c_p=\sum_{k=0}^p
\frac{(p!)^2}{(k!)^2}(2k)!4^{p-k}$,
and $\frac{\omega_I^k}{k!}$
defined above can be compared:
\[
c_p\Theta^p \preceq\frac{\omega_I^{2p}}{(2p)!}
\]
because the faces of $c_p\Theta^p$ are
quaternionic subspaces in $TM$, and the faces of
$\frac{\omega_I^{2p}}{(2p)!}$ are complex subspaces
(\ref{_triana_calibra_Theorem_}).

Let $\rho$ be a precalibration on a complex manifold
(\ref{_precalibra_Definition_}), and $\rho^{p,p}$ be
its $(p,p)$-part. We show that a plane
$V\subset TM$ is a face of $\rho^{p,p}$ if and only for $\zeta(V)$ is a face of $\rho$ for all $\zeta\in U(1)$,
for the standard $U(1)$-action on $TM$ (\ref{_avera_cali_Theorem_}).

Applying this result to the special Lagrangian calibration on $(M,J)$ defined in
\cite{_Harvey_Lawson:Calibrated_} (see also \cite{_McLean:SpLag_}), we obtain the form $\Psi_{n}$, $n=\dim_{\Bbb
H}M$, which calibrates complex analytic Lagrangian subvarieties on
$(M,I)$ (these subvarieties are known to be
special Lagrangian on $(M,J)$; see e.g. \cite{Hit1}).
This argument is not hard to generalize to arbitrary dimension.

\hfill

In most cases listed in \cite{_Harvey_Lawson:Calibrated_} and elsewhere, a calibration form is parallel with
respect to the Levi-Civita connection. An interesting side effect of our construction of holomorphic Lagrangian
calibrations is an appearance of a family of calibrations which are not parallel, under any torsionless connection (\ref{_conne_not_exist_Claim_}). These
calibrations are associated with the so-called HKT structures in hypercomplex geometry. In physics the HKT manifolds appear as target manifolds with $N=(4,0)$ supersymmetric $\sigma$-models with Wess-Zumino term \cite{_Howe_Papado_}.
 
 We construct calibrations on a special class of hypercomplex manifolds with holonomy of its Obata connection in $SL(n, {\Bbb H})$, the commutator subgroup of $GL(n,{\Bbb H})$. Such manifolds are called {\bf $SL(n, {\Bbb H})$-manifolds}. For more
examples and an introduction to $SL(n, {\Bbb H})$-geometry, see Section \ref{_SL_n_H_Section_}.
For any $SL(n, {\Bbb H})$-manifold $M$, and an induced complex structure $I$, there is a holomorphic volume form
$\Phi\in \Lambda^{2n,0}(M,J)$, which is parallel with respect to the Obata connection
(\cite{_Verbitsky:canoni_},
\cite{_BDV:nilmanifolds_}). The space $V$ of parallel
holomorphic volume forms is 1-dimensional. A choice of an
auxiliary induced complex structure such that $I\circ J =
- J\circ I$ endows $V$ with a real structure and a
positive direction (Subsection
\ref{_posi_2,0-forms_Subsection_}). We choose $\Phi$
to be real and positive.
Denote by $\Pi^{n,n}_I$ the projection to
$(n,n)$-component with respect to the complex structure
$I$, such that $I\circ J = - J\circ I$. 

\hfill

In Section 6 we show that
$\Re(\Pi^{n,n}_I\Phi)$ is a calibration for any
quaternionic Hermitian metric $g$ for which $|\Phi|=2^n$
(\ref{_Lagra_calibra_on_SL_n_H_Theorem_}).
This calibration calibrates complex subvarieties of $Z\subset (M,I)$
which are Lagrangian with respect to the $(2,0)$-form $\Omega=\omega_J+\1\omega_K$, defined as in
\eqref{_2,0_Equation_}.

This calibration is defined for any quaternionic Hermitian
metric, subject to the condition $|\Phi|=1$ (and
there are always many).  When $(M, I, J, K, \Phi, g)$ is
an HKT manifold with $\Hol(M)\subset \SLnH$, more
calibrations can be defined.

We choose $\Phi$ to be positive, real $(2n,0)$-form on
$(M,J)$, and let $\Phi_n:=\Re\Pi^{n,n}_I(\Phi)$. In
\cite{_Verbitsky:balanced_} it was shown that the form
$\Phi_{n+k}:= \frac 1 {2^kk!}\Phi_n \wedge \omega_I^k$ is always closed
and positive (\ref{_V_main_Proposition_}).
In \ref{_coisotro_calibra_on_balanced_Theorem_}, we prove
that this form is a calibration, for a metric $g':=g\cdot
\left|\frac{\Phi_{n+k}}{2^n}\right|^{(2n+2k)^{-1}}$, conformally equivalent to
$g$. When $g$ is also balanced, $|\Phi|=const$,
the conformal weight
$\left|\frac{\Phi_{n+k}}{2^n}\right|^{(2n+2k)^{-1}}$ is constant
(\ref{_Lagra_calibra_on_SL_n_H_Theorem_}),
and $g'$ is also HKT, but otherwise $g'$ is not an HKT
metric. In either case, the calibration $\Phi_{n+k}$
is (generally speaking) not parallel with respect to
any connection on $M$ (\ref{_conne_not_exist_Claim_}).

We show that $\Phi_{n+k}$ calibrates complex subvarieties
of $(M,I)$ which are coisotropic with respect to the
(2,0)-form $\Omega=\omega_J+\1\omega_K$
(\ref{_coisotro_calibra_on_HKT_Theorem_}).
The situation with {\em isotropic} subvarieties
is completely different. Using the
examples from Section 5.6, we notice in Remark 6.5 that
complex isotropic submanifolds in this case do not have to be
calibrated by any form, since they could be homologous to
zero.


\section{Preliminaries}


\subsection{Calibrations in Riemannian geometry}


\hfill

We provide here the basic definitions of the theory of calibrations which we use in the paper. The standard reference for this material is \cite{_Harvey_Lawson:Calibrated_} and the reader may also consult \cite{_Joyce:Calibrated_} for recent progress and developments related to manifolds with restricted holonomy.

\hfill

\definition
 Let $W\subset V$ be a $p$-dimensional subspace in a Euclidean space,
and $\Vol(W)$ denote the Riemannian volume form of $W \subset V$, defined up to a sign. For any $p$-form $\eta
\in \Lambda^p V$, let {\bf comass} $\comass(\eta)$ be  the maximum of $\frac{\eta(v_1, v_2, ...,
v_p)}{|v_1||v_2|...|v_p|}$, for all $p$-tuples $(v_1, ..., v_p)$ of vectors in $V$ and {\bf face} be the set of
planes $W\subset V$ where $\frac\eta {\Vol(W)}=\comass(\eta)$.

\hfill

\definition\label{_precalibra_Definition_}
A {\bf precalibration} on a Riemannian manifold is a differential form with comass $\leq 1$ everywhere.

\hfill

\definition
A {\bf calibration} is a precalibration which is closed.

\hfill

\definition
Let $\eta$ be a $k$-dimensional precalibration on a Riemannian manifold, and $Z\subset M$ a $k$-dimensional
subvariety (we usually assume that the Hausdorff dimension of the set of singular points of $Z$ is $\leq k-2$,
because in this case a compactly supported differential form can be integrated over $Z$). We say that $Z$ is
{\bf calibrated by $\eta$} if at any smooth point $z\in Z$, the space $T_zZ$ is a face of the precalibration
$\eta$.

\hfill

\remark Clearly, for any precalibration $\eta$,
\begin{equation}\label{_calibra_mini_Equation_}
\Vol(Z) \geq \int_Z\eta,
\end{equation}
 where
$\Vol(Z)$ denotes the Riemannian volume of a compact $Z$, and the equality happens iff $Z$ is calibrated by $\eta$. If, in
addition, $\eta$ is closed, $\int_Z\eta$ is a cohomological invariant, and the inequality
\eqref{_calibra_mini_Equation_} implies that $Z$ minimizes the Riemannian volume in its homology class.


\subsection{Hyperk\"ahler manifolds and calibrations}

The following definitions are standard.

\hfill

\definition
A manifold $M$ is called {\bf hypercomplex} if $M$ is equipped with a triple of complex structures $I, J, K$,
satisfying the quaternionic relations $I\circ J =- J \circ I =K$. If, in addition, $M$ is equipped with a
Riemannian metric $g$ which is K\"ahler with respect to $I, J, K$,
$(M, I, J, K, g)$ is called {\bf  hyperk\"ahler}.
This is equivalent to $\nabla I=\nabla J = \nabla K =0$, where $\nabla$ is the Levi-Civita connection of $g$;
see \cite{_Besse:Einst_Manifo_}.

\hfill

\remark\label{_Obata_Remark_} It has been known since 1955 that any hypercomplex manifold admits a torsion-free
connection preserving $I, J$ and $K$, which is necessarily unique. This connection is called {\bf the Obata
connection}, after M. Obata, who discovered it in \cite{_Obata_}. Any almost complex structure which is
preserved by a torsion-free connection is necessarily integrable (this is an easy consequence of
Newlander-Nirenberg theorem). Therefore, for any $a, b, c \in \R$, with $a^2 + b ^2 + c^2=1$, the almost complex
structure $aI + b J + c K$ is in fact integrable. We denote by $(M,L)$ the manifold $M$ considered as a complex manifold
with the complex structure induced by $L=aI + b J + c K$.

\hfill

\definition
Such complex structures are called {\bf induced by quaternions}, and the corresponding family, parametrized by
$S^2$ -- {\bf the twistor family}, or {\bf the hypercomplex family}. This family is holomorphic, and its total
space (fibered over $\C P^1$) is called {\bf the twistor space of $M$}. It is a complex analytic space,
non-K\"ahler even in simplest cases (for $M$ a torus or a K3 surface).

\hfill

Hyperk\"ahler geometry has a long history and is already well established. For more details and background
definitions, please see \cite{_Besse:Einst_Manifo_,_Joyce:Calibrated_}. In algebraic geometry, the word {\em hyperk\"ahler} is
essentialy synonymous with ``holomorphic symplectic''. The reason is that any hyperk\"ahler manifold is equipped
with a complex-valued  form
$\Omega:=\omega_J+\1\omega_K$.%
\footnote{We always write $\omega_I, \omega_J, \omega_K$ for the corresponding K\"ahler forms.} This form has
Hodge type (2,0) on $(M,I)$ and is closed, hence holomorphically symplectic.

The converse follows from the Yau's proof of Calabi's conjecture: a holomorphically symplectic, K\"ahler
manifold admits a unique hyperk\"ahler metric in a given K\"ahler class (\cite{_Besse:Einst_Manifo_}). For survey of recent advances in hyperk\"ahler geometry see \cite{_Huybrechts1_,_Huybrechts2_}.

\hfill

Some of the main objects of this paper are holomorphic
Lagrangian, isotropic and coisotropic subvarieties of
$(M,I)$, where $(M, I, J, K, g)$ is hyperk\"ahler.

\hfill

\definition
A complex analytic subvariety $Z$ of a holomorphically symplectic manifold $(M, \Omega)$ is called {\bf
holomorphic Lagrangian} if $\Omega\restrict Z=0$, and $\dim_\C Z = \frac 1 2 \dim_\C M$, and {\bf isotropic} if
$\Omega\restrict Z=0$, and $\dim_\C Z < \frac 1 2 \dim_\C M$. It is called {\bf coisotropic} if $\Omega$ has
rank $\frac 1 2 \dim_\C M -\codim_\C Z$ on $TZ$ in all smooth points of $Z$, which is the minimal possible rank
for a $2n-p$-dimensional subspace in a $2n$-dimensional symplectic space.


\subsection{Calibrations in HKT-geometry}

Let $(M, I, J, K)$ be a hypercomplex manifold. Then the tangent bundle $TM$ is equipped with a natural
quaternionic action. In particular, the group $SU(2)$ of unitary quaternions acts on $TM$, in a canonical way. A
Riemannian metric on $M$ is called {\bf quaternionic Hermitian} if it is $SU(2)$-invariant. A hyperk\"ahler
metric is obviously quaternionic Hermitian, but the converse is manifestly false, as we shall explain presently.

With every quaternionic Hermitian metric $g$ we associate 2-forms $\omega_I:= g(I\cdot, \cdot), \omega_J:=
g(J\cdot, \cdot)$ and $\omega_K:= g(K\cdot, \cdot)$ which are clearly antisymmetric, because $g$ is
$SU(2)$-invariant. It is easy to check that
\begin{equation}\label{_2,0_Equation_}
\Omega:= \omega_J + \1 \omega_K
\end{equation}
is a (2,0)-form on $(M,I)$. This form is closed if and only if $(M,I,J,K,g)$ is hyperk\"ahler
(\cite{_Besse:Einst_Manifo_}).

For a weaker form of this condition, consider the $(1,0)$-part of the de Rham differential,
\[
\6:\; \Lambda^{p,q}(M,I)\arrow \Lambda^{p+1,q}(M).
\]
A quaternionic Hermitian hypercomplex manifold is called
{\bf HKT} (short for ``hyperk\"ahler with torsion'') if
$\6\Omega=0$.

The theory of HKT-manifolds is a rapidly developing
subfield of quaternionic geometry. Originally this notion
appeared in physics (\cite{_Howe_Papado_}), but
mathematicians found it very useful. For an early survey
of HKT-geometry, please see \cite{_Gra_Poon_}.

Another ingredient of an HKT calibration theory is the
notion of Obata connection (\ref{_Obata_Remark_}). Since
this connection preserves the quaternionic structure, its
holonomy $\Hol(M)$ lies in $GL(n, {\Bbb H})$. The
holonomy of the Obata connection is one of the most
important invariants of a hypercomplex manifold. Many
properties of $M$ can be related directly to its
holonomy group. In particular, the group $\Hol(M)$ is
compact if and only if $(M,I,J,K)$ admits a hyperk\"ahler metric.

There seems to be no holonomy characterization of HKT structures. In fact
the holonomy of Obata connection is rarely known explicitly, except on hyperk\"ahler manifolds, where it is equal to the Levi-Civita connection.
However the knowledge of holonomy is still quite useful for the study of HKT geometry.
For many examples of compact  hypercomplex manifolds, the
group $\Hol(M)\subset GL(n, {\Bbb H})$
is strictly smaller than $GL(n, {\Bbb H})$.
Only recently it was found that the group $SU(3)$
with the left-invariant hypercomplex structure
has $GL(n, {\Bbb H})$ as its holonomy group
(\cite{_Soldatenkov:SU(2)_}).

An important
subgroup inside $GL(n, {\Bbb H})$ is its commutator
$SL(n, {\Bbb H})$. This group can be defined as a group of
quaternionic matrices $A \subset \End({\Bbb H}^n)$
preserving a non-zero complex-valued form $\Phi\in \Lambda^{2n,0}_\C({\Bbb H}^n_I)$, where ${\Bbb H}^n_I$ is
${\Bbb H}^n$ considered as a $2n$-dimensional complex space, with the complex structure $I$ induced by
quaternions. The coefficient $\lambda :=\frac{A(\Phi)}{\Phi}$ is called {\bf the Moore determinant} of the
matrix $A$ (\cite{_Alesker:non-commu_},
\cite{_Alesker_Verbitsky_HKT_});
it is always a positive real number, with
$\lambda^4$ equal to the determinant of $A$, considered as an element of $GL(4n, \R)$. The group $SL(n, {\Bbb
H})$ is a group of quaternionic matrices with Moore determinant 1.



\section{$SL(n, {\Bbb H})$-manifolds}
\label{_SL_n_H_Section_}


\subsection{An introduction to $SL(n, {\Bbb H})$-geometry}

As Obata has shown (\cite{_Obata_}), a hypercomplex manifold $(M,I,J,K)$ admits a necessarily unique
torsion-free connection, preserving $I,J,K$. The converse is also true: if a manifold $M$ equipped with an
action of ${\Bbb H}$ admits a torsion-free connection preserving the quaternionic action, it is hypercomplex.
This implies that a hypercomplex structure on a manifold can be defined as a torsion-free connection with
holonomy in $GL(n, {\Bbb  H})$. This connection is called {\bf the Obata connection} on a hypercomplex manifold.

Connections with restricted holonomy are one of the central notions in Riemannian geometry, due to Berger's
classification of irreducible holonomy of Riemannian manifolds. However, a similar classification exists for
general torsion-free connections (\cite{_Merkulov_Sch:long_}). In the Merkulov-Schwachh\"ofer list, only three
subroups of $GL(n, {\Bbb H})$ occur. In addition to the compact group $Sp(n)$ (which defines hyperk\"ahler
geometry), also $GL(n, {\Bbb H})$ and its commutator $SL(n, {\Bbb H})$ appear, corresponding to hypercomplex
manifolds and hypercomplex manifolds with trivial determinant bundle, respectively. Both of these geometries are
interesting, rich in structure and examples, and deserve detailed study.

It is easy to see that $(M,I)$ has holomorphically trivial canonical bundle, for any $SL(n, {\Bbb H})$-manifold
$(M,I, J, K)$ (\cite{_Verbitsky:canoni_}). For a hypercomplex manifold with trivial canonical bundle admitting
an HKT metric, a version of Hodge theory was constructed (\cite{_Verbitsky:HKT_}). Using this result, it was
shown that a compact hypercomplex manifold with trivial canonical bundle has holonomy in $SL(n,{\Bbb H})$, if it
admits an HKT-structure (\cite{_Verbitsky:canoni_}).

In \cite{_BDV:nilmanifolds_}, it was shown that holonomy of all hypercomplex nilmanifolds lies in $SL(n, {\Bbb
H})$. Many  working examples of hypercomplex manifolds are in fact nilmanifolds, and by this result they all
belong to the class of $SL(n, {\Bbb H})$-manifolds.

The $SL(n, {\Bbb H})$-manifolds were studied in \cite{_AV:Calabi_} and \cite{_Verbitsky:skoda.tex_}, because on
such manifolds the quaternionic Dolbeault complex is identified with a part of de Rham complex
(\ref{_V_main_Proposition_}).
 Under this identification, ${\Bbb H}$-positive forms become positive
in the usual sense, and $\6$, $\6_J$-closed or exact forms become $\6, \bar\6$-closed or exact (see Section 3.1). This
linear-algebraic identification is especially useful in the study of
the quaternionic Monge-Amp\`ere equation
(\cite{_AV:Calabi_}).

\subsection{Balanced HKT-manifolds}

The following lemma is contained in
\cite{_BDV:nilmanifolds_} (Theorem 3.2; see also
\cite{_Verbitsky:balanced_}, Lemma 4.3). Recall that
the map $\eta \arrow J(\bar \eta)$
defines a real structure on $\Lambda^{2p,0}(M,I)$.
A $(p,0)$-form $\eta$ is called {\bf ${\Bbb H}$-real}
if $J(\eta)=\bar\eta$.

\hfill

\lemma\label{_real_holo_parallel_Lemma_}
Let $(M,I,J,K)$
be a hypercomplex manifold, and $\eta$ a top degree
$(2n, 0)$-form, which is ${\Bbb H}$-real and
holomorphic. Then $\eta$ is Obata-parallel.

\endproof

\hfill

\definition
Let $(M, I, g)$ be a complex Hermitian manifold, $\dim_\C M =n$, and $\omega\in \Lambda^{1,1}(M)$ its Hermitian
form. One says that $M$ is {\bf balanced} if $d (\omega^{n-1})=0$.

\hfill

\remark It is easy to see that $d (\omega^{m})=0$ for $1 \leq m \leq n-2$ implies that $\omega$ is K\"ahler; the
balancedness makes sense as the only non-trivial condition of form $d (\omega^{m})=0$ which is not equivalent to
the K\"ahler property.

\hfill

\theorem \label{_balanced_CY_equiv_Theorem_}
Let $(M,I,J,K, \Omega)$ be an HKT-manifold as in Section 2.3, $\dim_{\Bbb H}M =n$.
If $\bar\6$ is the standard Dolbeault operator on $(M,I)$, then the following conditions are equivalent.
\begin{description}
\item[(i)] $\bar\6(\Omega^n)=0$ \item[(ii)] $\nabla(\Omega^n)=0$, where $\nabla$ is the Obata connection
\item[(iii)] The manifold $(M,I)$ with the induced quaternionic Hermitian metric is balanced as a Hermitian
manifold:
\[
d(\omega_I^{2n-1})=0.
\]
\end{description}
{\bf Proof:} \cite{_Verbitsky:balanced_}, Theorem 4.8. \endproof

\hfill

\remark A balanced HKT-manifold has holonomy in $SL(n, {\Bbb
  H})$.
This statement follows immediately from the implication (iii) $\Rightarrow$ (ii) of
\ref{_balanced_CY_equiv_Theorem_}. However the balanced HKT condition is a little stronger. It is shown in \cite{_Ivanov_Petkov_} that an HKT manifold has (restricted) holonomy of the Obata connection in $SL(n, {\Bbb H})$ if and only if it is (locally) conformally balanced.

\hfill

\remark The condition $\nabla(\Omega^n)=0$ is independent from the choice of a basis $I,J,K$, $IJ = -JI=K$ of
${\Bbb H}$. Indeed, suppose that $g\in SU(n)$, and $(I_1, J_1, K_1)= (g(I), g(J), g(K))$ is a new basis in
${\Bbb H}$. The corresponding HKT-form $\Omega_1 = \omega_{J_1}+ \1 \omega_{K_1}$ can be expressed as $\Omega_1
= g(\Omega)$, hence
\[
\nabla(\Omega_1^n)= \nabla(g(\Omega_1^n)) = g(\nabla(\Omega^n))=0.
\]
Therefore, \ref{_balanced_CY_equiv_Theorem_} leads to the following corollary.

\hfill

\corollary Let $(M,I,J,K, \Omega)$ be an HKT-manifold, such that the corresponding complex Hermitian manifold
$(M,I)$ is balanced. Then $(M, I_1)$ is balanced for any complex structrure $I_1$ induced by the quaternions.
Moreover, $(M,I,J,K, \Omega)$ is an $SL({\Bbb H}, n)$-manifold.
\endproof


\section{Differential forms on hypercomplex manifolds}


In this section, we give an introduction to the linear algebraic structures on the de Rham algebra of a
hypercomplex manifold. We follow \cite{_Verbitsky:skoda.tex_} and \cite{_Verbitsky:balanced_}.

\subsection{The quaternionic
  Dolbeault complex}

It is well-known that any irreducible representation of $SU(2)$ over $\C$ can be obtained as a symmetric power
$S^i(V_1)$, where $V_1$ is a fundamental 2-dimensional representation. We say that a representation $W$ {\bf has
weight $i$} if it is isomorphic to $S^i(V_1)$. A representation is said to be {\bf pure of weight $i$} if all
its irreducible components have weight $i$.

\hfill

\remark\label{_weight_multi_Remark_} The Clebsch-Gordan formula (see \cite{_Humphreys_}) claims that the weight
is {\em multiplicative}, in the following sense: if $i\leq j$, then
\[
V_i\otimes V_j = \bigoplus_{k=0}^i V_{i+j-2k},
\]
where $V_i=S^i(V_1)$ denotes the irreducible representation of weight $i$.

\hfill

Let $M$ be a hypercomplex  manifold, $\dim_{\Bbb H}M=n$. There is a natural multiplicative action of
$SU(2)\subset {\Bbb H}^*$ on $\Lambda^*(M)$, associated with the hypercomplex structure.

\hfill

Let $V^i\subset \Lambda^i(M)$ be a maximal $SU(2)$-invariant subspace of weight $<i$. The space $V^i$ is well
defined, because it is a sum of all irreducible representations $W\subset \Lambda^i(M)$ of weight $<i$. Since
the weight is multiplicative (\ref{_weight_multi_Remark_}), $V^*= \bigoplus_i V^i$ is an ideal in
$\Lambda^*(M)$.

It is easy to see that the de Rham differential $d$ increases the weight by 1 at most. Therefore, $dV^i\subset
V^{i+1}$, and $V^*\subset \Lambda^*(M)$ is a differential ideal in the de Rham DG-algebra $(\Lambda^*(M), d)$.

\hfill

\definition\label{_qD_Definition_}
Denote by $(\Lambda^*_+(M), d_+)$ the quotient algebra $\Lambda^*(M)/V^*$. It is called {\bf the quaternionic
Dolbeault algebra of
  $M$}, or {\bf the quaternionic Dolbeault complex}
(qD-algebra or qD-complex for short).

\hfill

\remark The complex $(\Lambda^*_+(M), d_+)$ was
constructed earlier by Capria and Salamon
(\cite{_Capria-Salamon_}) in a different (and more
general) situation, and much studied since then.

\hfill

The Hodge bigrading is compatible with the weight decomposition of $\Lambda^*(M)$, and gives a Hodge
decomposition of $\Lambda^*_+(M)$ (\cite{_Verbitsky:HKT_}):
\[
\Lambda^i_+(M) = \bigoplus_{p+q=i}\Lambda^{p,q}_{+,I}(M).
\]
The spaces $\Lambda^{p,q}_{+,I}(M)$ are the weight spaces for a particular choice of a Cartan subalgebra in
$\goth{su}(2)$. The $\goth{su}(2)$-action induces an isomorphism of the weight spaces within an irreducible
representation. This gives the following result.

\hfill

\proposition \label{_qD_decompo_expli_Proposition_} Let $(M,I,J,K)$ be a hypercomplex manifold and
\[
\Lambda^i_+(M) = \bigoplus_{p+q=i}\Lambda^{p,q}_{+,I}(M)
\]
the Hodge decomposition of qD-complex defined above. Then there is a natural isomorphism
\begin{equation}\label{_qD_decompo_Equation_}
\Lambda^{p,q}_{+,I}(M)\cong \Lambda^{p+q,0}(M,I).
\end{equation}

{\bf Proof:} See \cite{_Verbitsky:HKT_}. \endproof

\hfill

This isomorphism is compatible with a natural algebraic structure on \[ \bigoplus_{p+q=i}\Lambda^{p+q,0}(M,I),\]
and with the Dolbeault differentials, in the following way.

\hfill

Let $(M,I,J,K)$ be a hypercomplex manifold. We extend \[ J:\; \Lambda^1(M) \arrow \Lambda^1(M)\] to
$\Lambda^*(M)$ by multiplicativity. Recall that
\[ J(\Lambda^{p,q}(M,I))=\Lambda^{q,p}(M,I), \]
because $I$ and $J$ anticommute on $\Lambda^1(M)$. Denote by
\[ \6_J:\;  \Lambda^{p,q}(M,I)\arrow \Lambda^{p+1,q}(M,I)
\]
the operator $J\circ \bar\6 \circ J$, where $\bar\6:\;  \Lambda^{p,q}(M,I)\arrow \Lambda^{p,q+1}(M,I)$ is the
standard Dolbeault operator on $(M,I)$, that is, the $(0,1)$-part of the de Rham differential. Since
$\bar\6^2=0$, we have $\6_J^2=0$. In \cite{_Verbitsky:HKT_} it was shown that $\6$ and $\6_J$ anticommute:
\begin{equation}\label{_commute_6_J_6_Equation_}
\{\6_J, \6 \}=0.
\end{equation}

Consider the quaternionic Dolbeault complex $(\Lambda^*_+(M), d_+)$ constructed in \ref{_qD_Definition_}. Using
the Hodge bigrading, we can decompose this complex, obtaining a bicomplex
\[
\Lambda^{*, *}_{+,I}(M) \xlongrightarrow{d^{1,0}_{+,I}, d^{0,1}_{+,I}} \Lambda^{*, *}_{+,I}(M)
\]
where $d^{1,0}_{+,I}$,  $d^{0,1}_{+,I}$ are the Hodge components of the quaternionic Dolbeault differential
$d_+$, taken with respect to $I$.

\hfill

\theorem\label{_bico_ide_Theorem_} Under the multiplicative isomorphism
\[
\Lambda^{p,q}_{+,I}(M)\cong \Lambda^{p+q,0}(M,I)
\]
constructed in \ref{_qD_decompo_expli_Proposition_}, $d^{1,0}_+$ corresponds to $\6$ and $d^{0,1}_+$ to $\6_J$:

\begin{equation}\label{_bicomple_XY_Equation}
\begin{minipage}[m]{0.85\linewidth}
{\tiny $ \xymatrix @C+1mm @R+10mm@!0  {
  && \Lambda^0_+(M) \ar[dl]^{d^{0,1}_+} \ar[dr]^{d^{1,0}_+}
   && && && \Lambda^{0,0}_I(M) \ar[dl]^{\6} \ar[dr]^{ \6_J}
   &&  \\
 & \Lambda^{1,0}_+(M) \ar[dl]^{d^{0,1}_+} \ar[dr]^{d^{1,0}_+} &
 & \Lambda^{0,1}_+(M) \ar[dl]^{d^{0,1}_+} \ar[dr]^{d^{1,0}_+}&&
\text{\large $\cong$} &
 &\Lambda^{1,0}_I(M)\ar[dl]^{ \6} \ar[dr]^{ \6_J}&  &
 \Lambda^{1,0}_I(M)\ar[dl]^{ \6} \ar[dr]^{ \6_J}&\\
 \Lambda^{2,0}_+(M) && \Lambda^{1,1}_+(M)
   && \Lambda^{0,2}_+(M)& \ \ \ \ \ \ & \Lambda^{2,0}_I(M)& &
\Lambda^{2,0}_I(M) & &\Lambda^{2,0}_I(M) \\
} $ }
\end{minipage}
\end{equation}
Moreover, under this isomorphism,
the form $\omega_I\in \Lambda^{1,1}_{+,I}(M)$ corresponds to
$\Omega\in\Lambda^{2,0}_I(M)$.

{\bf Proof:} See \cite{_Verbitsky:HKT_} or \cite{_Verbitsky:qD_}.
 \endproof

\subsection{Positive $(2,0)$-forms on hypercomplex
  manifolds}
\label{_posi_2,0-forms_Subsection_}

The notion of positive $(2p,0)$-forms on hypercomplex manifolds (sometimes called q-positive, or ${\Bbb
H}$-positive) was developed in \cite{_Alesker_Verbitsky_HKT_} (see also \cite{_AV:Calabi_} and
\cite{_Verbitsky:skoda.tex_}).

\hfill

Let $\eta\in \Lambda^{p,q}_I(M)$ be a differential form. Since $I$ and $J$ anticommute, $J(\eta)$ lies in
$\Lambda^{q,p}_I(M)$. Clearly, $J^2\restrict {\Lambda^{p,q}_I(M)}=(-1)^{p+q}$. For $p+q$ even, $J\restrict
{\Lambda^{p,q}_I(M)}$ is an anticomplex involution, that is, a real structure on $\Lambda^{p,q}_I(M)$. A form
$\eta \in \Lambda^{2p,0}_I(M)$ is called {\bf real} if $J(\bar\eta)=\eta$.

For a real $(2,0)$-form $\eta$,
\[ 
   \eta\left(x, J(\bar x))\right)=
   \bar \eta\left(J(x), J^2 (\bar x)\right)=
 \bar \eta\left(\bar x, J(x)\right),
\] 
for any $x \in T^{1,0}_I(M)$. {}{}From the definition of a real form, we obtain that the scalar $\eta\left(x,
J(\bar x)\right)$ is always real.

\hfill

\definition
A real $(2,0)$-form $\eta$ on a hypercomplex manifold is called {\bf positive} if $\eta\left(x, J(\bar
x)\right)\geq 0$ for any $x \in T^{1,0}_I(M)$, and {\bf strictly positive} if this inequality is strict, for all
$x\neq 0$.

\hfill

An HKT-form $\Omega\in \Lambda^{2,0}_I(M)$ of any
HKT-structure is strictly positive. Moreover,
HKT-structures on a hypercomplex manifold are in one-to-one
correspondence with $\6$-closed, strictly positive
$(2,0)$-forms.

The analogy between K\"ahler forms and HKT-forms can be pushed further; it turns out that any HKT-form
$\Omega\in \Lambda^{2,0}_I(M)$ has a local potential $\phi\in C^\infty(M)$, in such a way that
$\6\6_J\phi=\Omega$ (\cite{_Alesker_Verbitsky_HKT_}). Here $\6\6_J$ is a composition of $\6$ and $\6_J$ defined
on quaternionic Dolbeault complex as above (these operators anticommute).

\subsection[The map ${\cal V}_{p,q}:\;
  \Lambda^{p+q,0}_I(M)\arrow\Lambda^{n+p, n+q}_I(M)$
on $SL(n, {\Bbb H})$-manifolds]{The map ${\cal V}_{p,q}:\;
  \Lambda^{p+q,0}_I(M)\arrow\Lambda^{n+p, n+q}_I(M)$\\
on $SL(n, {\Bbb H})$-manifolds} \label{_V_p,q_Subsection_}

Let $(M,I,J,K)$ be an $SL(n, {\Bbb H})$-manifold, $\dim_\R M =4n$, and
\[
  {\cal R}_{p,q}:\; \Lambda^{p+q,0}_I(M)\arrow \Lambda^{p,q}_{I,+}(M)
\]
the isomorphism induced by $\goth{su}(2)$-action as in \ref{_bico_ide_Theorem_}. Consider the projection
\begin{equation}\label{_proj_to_+_Equation_}
\Lambda^{p,q}_{I}(M)\arrow \Lambda^{p,q}_{I,+}(M),
\end{equation}
and let $R:\; \Lambda^{p,q}_{I}(M)\arrow\Lambda^{p+q,0}_I(M)$ denote the composition of
\eqref{_proj_to_+_Equation_} and ${\cal R}_{p,q}^{-1}$.

Let $\Phi_I$ be a nowhere degenerate holomorphic section of $\Lambda^{2n,0}_I(M)$. Assume that $\Phi_I$ is real,
that is, $J(\Phi_I)=\bar\Phi_I$, and positive. Existence of such a form is equivalent to $\Hol(M) \subset SL(n,
{\Bbb H})$ (\ref{_real_holo_parallel_Lemma_}). It is often convenient to define $SL(n, {\Bbb H})$-structure by
fixing the quaternionic action and the holomorphic form $\Phi_I$.

\hfill

Define the map
\[ {\cal V}_{p,q}:\;
  \Lambda^{p+q,0}_I(M)\arrow\Lambda^{n+p, n+q}_I(M)
\]
by the relation
\begin{equation}\label{_V_p,q_via_test_form_Equation_}
{\cal V}_{p,q}(\eta) \wedge \alpha = \eta \wedge R(\alpha)\wedge \bar\Phi_I,
\end{equation}
for any test form $\alpha \in \Lambda^{n-p, n-q}_I(M)$.

\hfill

The map ${\cal V}_{p,p}$ is especially remarkable, because it maps closed, positive $(2p,0)$-forms to closed,
positive $(n+p, n+p)$-forms, as the following proposition implies.

\hfill

\proposition\label{_V_main_Proposition_}
Let $(M,I,J,K, \Phi_I)$ be an $SL(n, {\Bbb H})$-manifold, and
\[ {\cal V}_{p,q}:\;
  \Lambda^{p+q,0}_I(M)\arrow\Lambda^{n+p, n+q}_I(M)
\]
 the map defined above.
Then
\begin{description}
\item[(i)] ${\cal V}_{p,q}(\eta)= {\cal R}_{p,q}(\eta)
  \wedge {\cal V}_{0,0}(1)$.

\item[(ii)]  The map ${\cal
V}_{p,q}$ is injective, for
  all $p$, $q$.

\item[(iii)] $(\1)^{(n-p)^2}{\cal V}_{p,p}(\eta)$ is real if and
  only $\eta\in\Lambda^{2p,0}_I(M)$ is real,
and positive if and only if $\eta$ is
positive.

\item[(iv)] ${\cal V}_{p,q}(\6\eta)= \6{\cal
V}_{p-1,q}(\eta)$, and ${\cal V}_{p,q}(\6_J\eta)=
  \bar\6{\cal  V}_{p,q-1}(\eta)$.

\item[(v)] ${\cal V}_{0,0}(1)= \lambda {\cal
  R}_{n,n}(\Phi_I)$, where $\lambda$ is a positive rational number,
depending only on the dimension $n$.
\end{description}

{\bf Proof:} See \cite{_Verbitsky:skoda.tex_}, Proposition 4.2, or \cite{_AV:Calabi_}, Theorem 3.6. \endproof

\hfill

\remark\label{_V_Omega^k_Remark_}
For the purposes of the present paper, we are interested
in \ref{_V_main_Proposition_} for the case
$\eta=\Omega^k$, where $\Omega$ is an HKT-form.
In this case, ${\cal R}_{p,p}(\Omega^k)$ is a projection
of $\omega_I^k$ to the component of maximal weight (see
\ref{_Pi_+_of_omega^n+p_clo_Proposition_}  below).
Now, ${\cal V}_{p,q}(\Omega^k)= {\cal R}_{p,q}(\Omega^k)
  \wedge {\cal V}_{0,0}(1)$, as follows from
  \ref{_V_main_Proposition_} (i). However,
${\cal V}_{0,0}(1)$ has weight $2n$, by
\ref{_V_main_Proposition_} (v), and
 $\omega_I^k$ has weight $\leq 2k$,
hence their product is of weight $\geq 2n-2k$.
Since this product is $(2n-2k)$-form, it is
pure of weight $(2n-2k)$, and components of
$\omega_I^k$ of weight $<2k$ do not contribute to
the product $\omega^k_I\wedge {\cal V}_{0,0}(1)$.
We obtain that the closed, positive form
${\cal V}_{k,k}(\Omega^k)$ is proportional
to $\omega^k_I\wedge {\cal V}_{0,0}(1)$, with positive
coefficient.

\subsection{Algebra generated by $\omega_I$, $\omega_J$,
  $\omega_K$}
\label{_alge_gene_omega_I_etc_Subsection_}

Let $(M,I,J,K, g)$ be a quaternionic Hermitian
manifold. Consider the algebra $A^*= \oplus A^{2i}$
generated by $\omega_I$, $\omega_J$, and $\omega_K$. In
\cite{_Verbitsky:trianalyt_}, this algebra was computed
explicitly. It was shown that, up to the middle degree,
$A^*$ is a symmetric algebra with generators $\omega_I$,
$\omega_J$, $\omega_K$. The algebra $A^*$ has Hodge
bigrading $A^k = \bigoplus\limits_{p+q=k}A^{p,q}$. {}From
the Clebsch-Gordan formula, we obtain that $A^{2i}_+:=
\Lambda^{2i}_+(M)\cap A^{2i}$, for $i\leq n$, is an
orthogonal complement to $Q(A^{2i-4})$, where $Q(\eta) =
\eta \wedge (\omega_I^2 + \omega_J^2+\omega_K^2)$.
Moreover, $A^{2i}_+$ is irreducible as a representation of
$SU(2)$. Therefore, the space $A^{p,p}_+= \ker
Q^*\restrict {A^{p,p}}$ is 1-dimensional.
This argument also implies that
the form ${\cal V}_{0,0}(1)$ is proportional
to $\Phi_J|^{n,n}_I$, where
$\Phi_J$ is a holomorphic volume form
on $(M,J)$, obtained as a top power
of the appropriate holomorphic symplectic
form, and $\Phi_J|^{n,n}_I$ its $(n,n)$-part,
taken with respect to $I$.

\hfill

\proposition\label{_Pi_+_of_omega^n+p_clo_Proposition_}
Let  $(M,I,J,K, \Phi_I)$ be an $SL(n, {\Bbb H})$-manifold,
equipped with an HKT-structure $\Omega$. Assume that $\Omega^n = \Phi_I$. Let
\[ \Pi_+:\; \Lambda^{n+k,n+k}_I(M)\arrow \Lambda^{n+k,n+k}_{I,+}(M)\]
be the projection to the component of maximal weight with respect to the $SU(2)$-action. Then
$\Xi_k:=\Pi_+(\omega_I^{n+k,n+k})$ is a closed, weakly
positive $(n+k,n+k)$-form, which is proportional to
$\omega^k_I\wedge \Phi_J|^{n,n}_I$ and to
$\omega^k_I\wedge {\cal V}_{0,0}(1)$.

\hfill

{\bf Proof:}
The form $\omega^k_I\wedge \Phi_J|^{n,n}_I$ is
proportional to $\omega^k_I\wedge {\cal V}_{0,0}(1)$
as indicated above.
Consider the algebra $A^*= \oplus A^{2i}$ generated by
$\omega_I$, $\omega_J$, and $\omega_K$. The map $R^{p,q}$
is induced by the $SU(2)$-action, hence it maps $A^{*,*}$
to itself. Since ${\cal V}_{p,q}(\eta)= {\cal
  R}_{p,q}(\eta) \wedge {\cal V}_{0,0}(1)$, and ${\cal
  V}_{0,0}(1)$ is proportional to
${\cal R}_{n,n}(\Phi_I)\in A^*$, we obtain
\[ {\cal V}_{p,q} (A^{p+q,0}) \subset A^{n+p,n+q}.
\]
Since ${\cal V}_{p,p}(\Omega^p)\subset A^{n+p,n+p}_+$,
the 1-dimensional space $A^{n+p,n+p}_+$ is generated by ${\cal
V}_{p,p}(\Omega^p)$. This form is closed and positive by
\ref{_V_main_Proposition_}. Therefore, the projection of
$\omega_I^{n+p}$ to $A^{n+p,n+p}_+$ is closed and positive
(see \ref{_V_Omega^k_Remark_}).
\endproof


\section{Calibrations on hyperk\"ahler manifolds}


\subsection{Hodge decomposition and $U(1)$-action}

Let $I$ be a complex structure on a vector space $V$ and $\rho:\; U(1) \arrow \End(V)$ a real
$U(1)$-representation given by $\rho(t)(X)=(\cos t + \sin t I)X$. This is extended by multiplicativity to a
representation in the tensor powers of $V$ with $\rho(t)(\alpha)(X)=\alpha(\rho(t)X)$ for a 1-form $\alpha$. In
the usual fashion, we define the weight decomposition
associated with
this $U(1)$-action: the tensor $z$ has weight $p$
if $\rho(t) z = (\cos pt)z + \1 (\sin pt) z$.
We need also the definition of average over $U(1)$ of $Y$:
\[
    \Av_{\rho}Y=\frac{1}{2\pi}\int_0^{2\pi}\rho(t)Y dt.
\]
Note that $\rho(t)Y=Y$ for every $t$ implies $IY=Y$ for any
tensor $Y$ and that $I\Av_{\rho}Y=\Av_{\rho}Y$.

\hfill

\lemma\label{_averaging_calibrations_Lemma_}
Let $\rho$ be a U(1)-action on $W$, and $W = \bigoplus W^i$ the
corresponding weight decomposition. Then the projection to
$W^0$ along the sum of other $W^i$, $i\neq 0$,
coincides with taking the average over $U(1)$.

{\bf Proof:} For each $\eta\in W^i$, $i\neq 0$, one has
$\int_{U(1)} \rho(t) \eta dt=0$, because $\int_0^{2\pi}\cos(t)
dt=0$. \endproof

\hfill

\theorem\label{_avera_cali_Theorem_}
Let $\eta$ be a $2p$-form on a complex vector space $W$, with
$\comass(\eta)\leq 1$, and $\eta^{p,p}=\Av_{\rho}\eta$ be the $(p,p)$-part of $\eta$. Then $\comass
(\eta^{p,p})\leq 1$. Moreover, a $2p$-dimensional plane $V$ is a face of $\eta^{p,p}$ if and only if
$\rho(t)(V)$ is a face of $\eta$ for all $t\in \R$.

\hfill

{\bf Proof:}  For any decomposable $2p$-vector $\xi$,
its image $\rho(t)(\xi)$ is again decomposable for any $t$ and
$|\rho(t)(\xi)|=|\xi|$. Then
\[
\eta^{(p,p)}(\xi)= (\Av_\rho(\eta))(\xi)= \frac{1}{2\pi}\int_0^{2\pi}\eta(\rho(t)(\xi))\leq 1
\]
since $\eta(\rho(t)\xi)\leq 1$ for every $t$. The equality holds iff $\eta(\rho(t)\xi) = 1$ for every $t$.
\endproof

\subsection{An $SU(2)$-invariant calibration}

The most obvious example of a calibration on a hyperk\"ahler manifold is provided by the following theorem (see \cite{Berger} for similar statement about a quaternionic Wirtinger's inequality).

\hfill

\theorem\label{_triana_calibra_Theorem_}
Let $(M,I,J,K,g)$ be a hyperk\"ahler manifold, $\omega_I, \omega_J,
\omega_K$ the corresponding symplectic forms, and
$\Theta_p:= \frac{(\omega_I^2+ \omega_J^2
  +\omega_K^2)^p}{c_p}$ the standard $SU(2)$-invariant
$4p$-form normalized by $c_p=\sum_{k=1}^p
\frac{(p!)^2}{(k!)^2}(2k)!4^{p-k}$.
Then $\Theta_p$ is a calibration, and its faces are
$p$-dimensional quaternionic subspaces of $TM$. Moreover,
the form $\Xi_p:=\frac{(\omega_J^2+
  \omega_K^2)^p}{(p!)^24^p}$
is also a calibration, with the same faces.

\hfill

{\bf Proof:} Consider the form
$\widetilde\Xi_p:= \frac{\omega_J^{2p}}{(2p)!}$. By
\ref{_averaging_calibrations_Lemma_}, ${(\tilde
  \Xi_p)}^{2p,2p}_I=\Xi_p$, where $(\cdot)^{2p,2p}_I$ is
an operation of taking $(2p,2p)$-part under the complex
structure $I$. Indeed, $\omega_J^{2p}
=\frac{(\Omega+ \bar  \Omega)^{2p}}{4^p}$,
where $\Omega$ is the standard $(2,0)$ form on $(M,I)$.
Then the $(2p,2p)$-part of $\omega_J^{2p}$ is equal to
\[ \frac{(2p)!}{(p!)^2}\frac{\Omega^p \wedge
  \bar\Omega^p}{4^p}= \frac{(2p)!(\omega_J^2 +
  \omega_K^2)^p}{(p!)^24^p}.
\]
 By \ref{_avera_cali_Theorem_}, a subspace $V\subset TM$
 is a face of $\Xi_p$ if and only if $\rho_I(t)(V)$ is a
 face of ${\widetilde\Xi_p}$ for all $t$, with $\rho_I(t)$
 the $U(1)$-action associated with $I$. The form ${\tilde
 \Xi_p}$ is a standard K\"ahler calibration associated
 with $J$; it follows from
 \cite{_Harvey_Lawson:Calibrated_} that
 $V\subset TM$ is a face of ${\widetilde\Xi_p}$ if and only
 if it is $J$-linear, that is, $\C$-linear with respect
 to the action of $\C$ induced by $J$. Since $\rho(t)(V)$
 is $J$-linear for all $t$, it remains $J$-linear if we
 act on $V$ by a group $G$ generated by $\rho_I$ and
 $\rho_J$, with $\rho_J$ a $U(1)$-action associated with
 $J$. Clearly, $G\cong SU(2)$ is the group of unitary
 quaternions acting on $\Lambda^*M$. Therefore, $V$ is a
 face of $\Xi_p$ if and only if $V$ is $g(J)$-linear, for
 all $g\in SU(2)$. This is equivalent to $V$ being a
 quaternionic subspace. Taking the average of $\Xi_p$ with
 respect to $SU(2)$ will not change its faces, because
 they are already $SU(2)$-invariant. Therefore,
 $\Av_{SU(2)}(\Xi_p)$ is a calibration with its faces
 quaternionic subspaces. Moreover it is
 $Sp(n)Sp(1)$-invariant $4p$-form, so it is proportional
 to $(\omega_I^2+\omega_J^2+\omega_K^2)^p$. Then, using
 \ref{_calibr_constants_} below, we obtain that that
 $\Av_{SU(2)}(\Xi_p) = \Theta_p$ by evaluating both forms
 on a fixed quaternionic subspace. \endproof

\hfill

\remark
Subvarieties calibrated by $\Theta_p$
are called {\bf trianalytic subvarieties}.
They were studied, at some length, in
\cite{_Verbitsky:trianalyt_} and \cite{_Verbitsky:desing_}.

\subsection{A holomorphic Lagrangian calibration}

\proposition\label{_complex_Lagrangian_subspaces_Proposition_}
Let $(V^{4p}, I,J,K, g)$ be a quaternionic Hermitian vector
space with fundamental forms $\omega_I,\omega_J,\omega_K$, and $\Psi\in \Lambda^{2p}(V)$ a $2p$-form  which is
the real part of
$\frac{1}{p!}(\omega_I-\1\omega_K)^p$
(it is a $(2p,0)$-form with respect to $J$).
Denote by $\Psi^{p,p}_I$ the $(p,p)$-part of $\Psi$ with respect to $I$. Then $\Psi^{p,p}_I$ has comass 1. Moreover,
a $2p$-dimensional subspace $W\subset V$ is calibrated by $\Psi^{p,p}_I$ if and only if $W$ is complex $I$-linear
and calibrated by $\Psi$.

\hfill

{\bf Proof:}
The real part of $\frac{1}{p!}(\omega_I-\1\omega_K)^p$
calibrates special Lagrangian subspaces taken with
respect to the symplectic form $\omega_J$ (see
\cite{_Harvey_Lawson:Calibrated_}). Therefore, any face of
$\frac{1}{p!}(\omega_I-\1\omega_K)^p$ is $\omega_J$-Lagrangian.
By \ref{_avera_cali_Theorem_}, a $2p$-dimensional plane $W$ is a face of $\Psi^{p,p}_I$ if and only
if $\rho(t)(W)$  is a face of $\Psi$ for all $t\in \R$. It follows by taking  $t=0$ that $W$ is
$\omega_J$-Lagrangian and by taking $t=\pi/2$ that $I(W)$ is $\omega_J$-Lagrangian too. But $I(W)$ is
$\omega_J$-Lagrangian iff $W$ is $\omega_K$-Lagrangian. By \cite{Hit1} (see also
\ref{_Lagra_then_linear_Remark_} below) $W$ determines an $I$-complex subspace.
\endproof

\hfill

\remark\label{_Lagra_then_linear_Remark_} Let $V$ be a quaternionic Hermitian space, $\dim_{\Bbb H} V =p$, and
$\xi \in \Lambda^{2p}V$ a decomposable $2p$-vector which is associated with a $2p$-dimensional subspace
$W\subset V$. Clearly, $W$ is Lagrangian with respect to $\omega_J$ if and only if $L_{\omega_J}\xi=0$ and
$\Lambda_{\omega_J}\xi=0$, where $L_{\omega_J}$, $\Lambda_{\omega_J}$ are the corresponding Hodge operators,
$L_{\omega_J}(\eta) :=\eta\wedge\omega_J$, and $\Lambda_{\omega_J}= *L_{\omega_J}*$ its Hermitian adjoint. If
$W$ is Lagrangian with respect to $J$ and $K$, one has
\begin{equation}\label{_comme_L_Lambda_zero_Equation_}
[L_{\omega_J},\Lambda_{\omega_K}]\xi=0.
\end{equation}
However, the commutator $[L_{\omega_J},\Lambda_{\omega_K}]$ acts on forms of type $(p,q)$ with respect to $I$ as
a multiplication by $(p-q) \1$ (see \cite{_so(5)_}). Then \eqref{_comme_L_Lambda_zero_Equation_} implies that
$\xi$ is of type $(p,p)$ with respect to $I$.

\hfill

\claim\label{_cal_V_and_Psi_Claim_}
Let $V$ be an $n$-dimensional quaternionic Hermitian space, and ${\cal
V}^{0,0}:\; \R \arrow \Lambda^{n,n}_I(V)$ be a map defined in Subsection \ref{_V_p,q_Subsection_} (in Subsection
\ref{_V_p,q_Subsection_} it was defined for $SL(n, {\Bbb H})$-manifolds, but the definition can be repeated for
quaternionic spaces word by word). Then ${\cal
  V}^{0,0}(1)=\Psi^{n,n}_I$,
where $\Psi^{n,n}_I$ is a form defined as in
\ref{_complex_Lagrangian_subspaces_Proposition_}.

\hfill

{\bf Proof:}
From \ref{_V_main_Proposition_} (v),
we know that ${\cal V}^{0,0}(1)$ and $\Psi^{n,n}_I$ are proportional
and we only have to calculate the coefficient
of proportionality. For this we use ${\cal
V}^{0,0}(1)\wedge \alpha = R(\alpha)\wedge
\overline{\Phi}_I$ for a particular choice of $\alpha$ as
\[ \alpha=\xi_1\wedge...\wedge \xi_n\wedge
\overline{\xi_{n+1}}\wedge...\overline{\xi_{2n}},
\] where
$\xi_i$ are orthogonal and of unit norm. Then
\[ R(\alpha)=\xi_1\wedge...\wedge
   J\overline{\xi_{n+1}}\wedge...J\overline{\xi_{2n}}.
\]
From here if ${\cal V}^{0,0}(1)=\lambda\Psi^{n,n}_I$, then
$\lambda=1$. \endproof

\hfill

Comparing \ref{_V_main_Proposition_} and
\ref{_cal_V_and_Psi_Claim_},  we find that the form
$\Psi^{n,n}_I$ is positive.

\hfill


\subsection{Isotropic and coisotropic calibrations}

A similar argument can be applied to other powers of $\Omega_J$.

\hfill

\proposition  \label{_iso_cali_Proposition_}
Consider an $n$-dimensional quaternionic Hermitian space $V$, and let $\Omega_J:=
\omega_I-\1\omega_K$ be the usual $(2,0)$-form on the
complex space $(V,J)$. When $p\leq n$ denote by $\Psi_p:=
\frac 1 {p!}\Re (\Omega_J^p)$, and let $\Psi^{p,p}_I$ be
its $(p,p)$-part taken with respect to $I$. Then
$\Psi^{p,p}_I$ has comass 1, and its faces are complex
isotropic subspaces of $(V,I)$

\hfill

{\bf Proof:} Let $W\subset V$ be a real $2p$-dimensional
subspace, and $W_1$ be the smallest complex subspace of $(V,J)$
containing $W$. Adding more vectors if necessary, we can
always assume that $\dim_\C W_1=2p$. Denote by $\xi$ the
decomposable $4p$-vector associated with $W_1$, and
$I(\xi)$ its image under the action of a quaternion
$I$. Then $\frac 1 {p!}\Omega_J^p$ is a $(2p,0)$-form on
$W_1$, proportional to the unit holomorphic
volume form $\Vol^{2p,0}(W_1)$ with a coefficient $\kappa$  which satisfies
\[
|\kappa|= \frac{(\xi, I(\xi))}{|\xi|^2}
\]
where $(,)$ is the induced scalar product.
By Cauchy-Schwarz inequality $|\xi|\leq 1$, where the equality
holds iff $I\xi=\xi$ or, equivalently, $W_1$ is quaternionic. Since
$\Vol^{2p,0}(W_1)$ has comass 1,
\[
  \comass\left(\frac 1 {p!}\Omega_J^p\right)\leq 1
\]
with equality if and only if $W_1$ is quaternionic. In the latter case, $W$ is a face of $\frac 1
{p!}\Omega_J^p$ if and only if $W$ is complex Lagrangian in $W_1$, as follows from
\ref{_complex_Lagrangian_subspaces_Proposition_}.
\endproof

\hfill

We provide also an expression of $\Psi^{p,p}$ as a polynomial of $\omega_I,\omega_J$ and $\omega_K$ for even
$p$.

\hfill

\proposition\label{_explicit_form_Proposition_}
Let $\Psi^{p,p}$ be the $(p,p)$ part with respect to $I$ of
$Re(\omega_I-\1\omega_K)^p$. Then
\[
\Psi^{p,p} = \sum_{k=0}^{q} \frac{(-1)^k}{4^k} \binom{p}{2k}\binom{2k}{k}
\omega_I^{p-2k}\wedge(\omega_K^2+\omega_J^2)^{k}
\]
where $q=\llcorner \frac p 2 \lrcorner$ is the greatest
integer not exceeding $\frac{p}{2}$.

\hfill

{\bf Proof:} First we notice that
\[
Re(\omega_I-\1\omega_K)^p = \sum_{k=0}^{\llcorner \frac p 2 \lrcorner} (-1)^{k}
\binom{p}{2k}\omega_I^{p-2k}\wedge\omega_K^{2k}.
\]
 Since $\omega_I^{p-2k}$ is of type $(p-2k,p-2k)$ with respect to $I$ we need to determine
the type of $\omega_K^{2k}$. To do this  we use the fact that
$\omega_K=\frac{1}{2}\Omega+\frac{1}{2}\overline{\Omega}$ is the decomposition of $\omega_K$ in $(2,0)+(0,2)$
parts with respect to $I$ where $\Omega=\omega_K+\1\omega_J$. Then
\[
\omega_K^{2k} = \frac{1}{4^k}\sum_{s=0}^{p-2k} \binom{2k}{s}\Omega^{s}\wedge\overline{\Omega}^{2k-s}
\]
and each term in the sum has degree $(2s,4k-s)$ with respect to $I$. So the only term which will contribute to
$\Psi^{p,p}$ above will be when $s=k$. Obviously the term is
$\frac{1}{4^{k}}\binom{2k}{k}\Omega^{k}\wedge\overline{\Omega}^{k}$. Then the proposition follows from
the fact that $\Omega\wedge\overline{\Omega} = \omega_K^2+\omega_J^2$.
\endproof

\hfill

Notice that one can take the imaginary part of
$\Omega_J^p$ instead of the real part. The resulting
calibrated subspaces are again complex isotropic. To
identify the complex coisotropic subspaces, however, one has
to be more careful.

\hfill

\proposition\label{_Phi^p,p_comass_Proposition_}
Consider an $n$-dimensional quaternionic Hermitian space $V$, and let $\Omega_J:=
\omega_I-\1\omega_K$ be the usual $(2,0)$-form on the
complex space $(V,J)$. Let $\Phi_p +\sqrt{-1}\Phi'_p:=
\frac{1}{2^pp!n!}(\Omega_J)^n\wedge\omega_I^p$, and
 $\Phi^{p,p}_I$ (resp. $\Phi'^{p,p}_I$) be the
corresponding  $(n+p,n+p)$-parts taken with
respect to $I$. Then $\Phi^{p,p}_I$
(resp. $\Phi'^{p,p}_I$) have comass 1
and their faces are complex coisotropic subspaces of $(V,I)$

\hfill

{\bf Proof:} First we notice that if a form $\alpha$ is
calibration, then its Hodge dual $*\alpha$ is again
calibration and its faces are
orthogonal complements to the faces of
$\alpha$. Then the form $*\Psi^{p,p}$ is a calibration
with faces $I$-complex $\Omega_J$-coisotropic subspaces. The same is true also if we consider the imaginary part
of $\Omega_J^p$ instead of $\Psi^p$. Then it remains to check that the complex form in the proposition is Hodge
dual to $\Omega_J^p$ up to a real constant. To this end we first notice that
$*\Omega_J^{n-p}=c_1\overline{\Omega_J^n}\wedge\Omega_J^p$
for a real positive constant $c_1$. Then
$\Phi^{p,p}+\sqrt{-1}\Phi'^{p,p}$ and
$\overline{\Omega_J^n}\wedge\Omega_J^p$ are both highest vectors in an
irreducible representation $A^{2n+2p}$ of $SU(2)$
(see Subsection \ref{_alge_gene_omega_I_etc_Subsection_}),
hence they are proportional
up to a complex constant. More explicitly we have:
\[
(\omega_I-\1\omega_K)^n\wedge(\omega_I-\1\omega_K)^p =
\Omega_J^n\wedge(2\omega_I-\Omega_J)^p
\]
\[
=(\Omega_J)^n\wedge\sum_{s=0}^p\binom{p}{s}(-\Omega_J)^s\wedge2^{p-s}\omega_I^{p-s}
\]
Since $\Omega_J^{n+s}=0$ for $s>0$ all terms in the sum above vanish except the first one.
Then
\[
(\omega_I-\1\omega_K)^n\wedge(\omega_I+\1\omega_K)^p =
(\omega_I-\1\omega_K)^n\wedge2^{p}\omega_I^{p}
\]

{}From here and \ref{_calibr_constants_} $ii)$ the proposition follows.
\endproof

\hfill

To calculate the comass of the forms above we need the following well-known preliminary Lemma:

\hfill

\lemma\label{_kaehler_calib_constants_}
If $(V^{2n},I,g)$ is an Hermitian vector space and $\omega$
is the fundamental 2-form, then for any subset
$X_1,...,X_{2k}$ of a given unitary basis
$(e_1,Ie_1,...,e_n,Ie_n)$ we have:

i) $\omega^k(X_1,....,X_{2k})=\pm k!$ if $span\{X_1,...,X_{2k}\}$ is complex and

ii) $\omega^k(X_1,....,X_{2k})=0$ otherwise.

\hfill

The proof of $i)$ is standard, while $ii)$ follows from
the definition of wedge product and the fact that
$\omega(X_i,X_j) \neq 0$ only if $IX_i=\pm X_j$.

\hfill

\lemma\label{_calibr_constants_}
Let $(V^{4n}, I, J, K, g)$ be a real vector space
with anti-\-com\-mu\-ting complex structures $I, J, K$ compatible
with the positive scalar product $g$. Denote by $\omega_I, \omega_J, \omega_k$ the fundamental 2-forms corresponding
to $I, J$ and $K$ respectively and $\Omega_I=\omega_J+\sqrt{-1}\omega_K$ be the standard $I$-complex symplectic
2-form.  Consider the form
$\Psi_I^n=Re(\omega_I+\sqrt{-1}\omega_J)^n|_I^{(n,n)}$,
where $|_I^{(n,n)}$ denotes the $(n,n)$ component with
respect to $I$. Then:
\begin{description}

\item[i)] $\Omega_I^n\wedge\overline{\Omega_I}^n=4^n(n!)^2\Vol$ for the volume form $\Vol$ on $V$.

\item[ii)] $(\omega_I^2+\omega_J^2+\omega_K^2)^n=c_n \Vol$ where $c_n=\sum_{k=0}^n \frac{(n!)^2}{(k!)^2}(2k)!4^{n-k}$

\item[iii)] $\omega_I^k\wedge\Psi^n = 2^k k!n!\Vol_{E_{n+k}}$, where $E_{n+k}$ is an
$(n+k)$-dimensional $I$-complex and $\omega_J$-coisotropic subspace.
\end{description}

\hfill

{\bf Proof:} Fix a quaternionic-Hermitian co-basis
\[ (e^1, Ie^1, Je^1, Ke^1, e^2,Ie^2,...,Ke^n)
\]
of $V^*$ so
that $\Vol =e^1\wedge...\wedge Ke^n$ and let $e_1, Ie_1,...,Ke_n$ be the dual basis of $V$. From the fact that $\Omega_I=\sum_i dz^i\wedge dw^i$ for coordinates
$dz_i=e^i+\sqrt{-1}Ie^i$ and $dw_i=Je^i+\sqrt{-1}Ke^i$, follows that $\Omega_I^n=n!dz^1\wedge dw^1...dz^n\wedge dw^n$.
Then to obtain $i)$ we notice that $dz_i\wedge d\overline{z_i}=-2\sqrt{-1}e^i\wedge Ie^i$ and $dw_i\wedge
d\overline{w_i}=-2\sqrt{-1}Je^i\wedge Ke^i$.

To prove $ii)$ we write
\[ (\omega_I^2+\omega_J^2+\omega_K^2)^n = (\omega_I^2
   +\Omega_I\wedge\overline{\Omega_I})^n = \sum_{k=0}^n
    \binom{n}{k}\omega_I^{2k}\wedge \Omega_I^{n-k}\wedge\overline{\Omega_I}^{n-k}.
\]
Then we consider the term
$\omega_I^{2k}\wedge \Omega_I^{n-k}\wedge
\overline{\Omega_I}^{n-k}$.
Let $s_i = e^i\wedge Ie^i + Je^i\wedge Ke^i$ and
$t_j=dz^j\wedge dw^j$, so $\omega_I=\sum s_i$ and
$\Omega_I=\sum t_j$. Then $s_i^3=s_it_i=t_i^2=0$ $s_i,
t_j$ commute and $s_i^2=2\Vol_i,
t_i\overline{t_i}=4\Vol_i$, where
$\Vol_i= e^i\wedge Ie^i\wedge Je^i\wedge Ke^i$.
Fix $n-k$ indexes $(i_{k+1},i_{k+2},...,i_n)$. Then notice
that in the product
$\omega_I^{2k}\wedge t_{i_{k+1}}t_{i_{k+2}}...t_{i_n}\wedge \overline{\Omega_I}^{n-k}$
the only non-vanishing terms are of the form
\[
s_{i_1}^2s_{i_2}^2...s_{i_k}^2t_{i_{k+1}}t_{i_{k+2}}...t_{i_n}
\overline{t_{i_{k+1}}t_{i_{k+2}}...t_{i_n}}
\]
for the complementary indexes $(i_1,...,i_k)$, such that
$(i_1,...,i_n)$ is a permutation of $(1,2...,n)$. Every
such product is equal to $2^k4^{n-k}\Vol$. Then we may
select $i_1=1,..,i_k=k, i_{k+1}=k+1,...,i_n=n$ and count
the number of terms corresponding to it; clearly, this
number does not depend on the choice of the permutation.
The number is the product of the coefficients in front of
$s_1^2...s_k^2$ $t_{k+1}...t_{n}$ and
$\overline{t_{k+1}}...\overline{t_n}$ in the expansions of
$(s_1+...s_k)^{2k}$ $(t_{k+1}+...+t_n)^{n-k}$ and
$(\overline{t_{k+1}}+...+\overline{t_n})^{n-k}$
respectively, which is
$\frac{(2k)!}{2^k}((n-k)!)^2$. Since there are
$\frac{n!}{k!(n-k)!}$ different choices for $n-k$ indexes,
we obtain
\[ (\omega_I^2+\omega_J^2+\omega_K^2)^n=\sum_{k=0}^n
   \frac{(n!)^2}{(k!)^2}(2k)!4^{n-k} \Vol
\] and $ii)$
follows.

To prove $iii)$ we notice that $Sp(n)$ acts transitively on
complex coisotropic subspaces of fixed dimension. Then
 we choose the coisotropic subspace $L$ spanned by
$e_1, Ie_1,...., e_n, Ie_n, Je_1 ,Ke_1,...., Je_k, Ke_k$. Let
$\Omega_K=\omega_I+\sqrt{-1}\omega_J$, $\alpha\in L$
a subspace spanned by $2n$ vectors and $\beta$ be a
subspace generated by  $2k$ vectors among
$e_1, Ie_1,...., e_n, Ie_n, Je_1 ,Ke_1,...., Je_k, Ke_k$. Since
$\Psi^n=Re(\Omega_K)|^{n,n}_I$, then
$\Psi^n|_{\alpha}=0$ if $\alpha$ contains a
quaternionic subspace or is not $I$-invariant.
Similarly, $\omega_I^k\bigg |_{\beta}=0$ if $\beta$ is
not $I$-invariant as follows from Lemma \ref{_kaehler_calib_constants_}.

From the calculations in
\cite{_Harvey_Lawson:Calibrated_} p. 88,
we have
\[ \Psi^n(e_1,Ie_1,...,e_n,Ie_n)=n!Re(dz_1\wedge....\wedge
    dw_n)(e_1,Ie_1,...e_n,Ie_n)=n!,
\]
and from \ref{_kaehler_calib_constants_} above,
$\omega_I^K(Je_1,Ke_1,...,Je_k,Ke_k)=k!$.
Then in the expression for
$\Psi^n\wedge\omega_I^k(e_1,Ie_1,....Je_k,Ke_k)$
the only non-vanishing summands are
 $\omega_I^K(Je_1,Ke_1,...,Je_k,Ke_k)$ and the terms where one or more pairs
$e_i, Ie_i$ are interchanged with $Je_i, Ke_i$. If we have exactly $s$ pairs interchanged, then there will be
$\binom l s$ terms each with value $n!k!$. So
\[
  \Psi^n\wedge\omega_I^k(e_1,Ie_1,....Je_k,Ke_k)=
  n!k!\left (1+k+\binom k 2+...+ \binom k k\right) =2^nn!k!,
\]
which proves the Lemma. Note that for $k=n$
the result fits with the case $i)$ and
the calculations in \ref{_Phi^p,p_comass_Proposition_}.
\endproof

\hfill

\subsection{Holomorphic Lagrangian calibrations of degree two}

The calibration 4-forms with constant coefficients in $\R^8$ were studied
systematically in \cite{_DHM:R^8_}. Also various 4-forms
which are calibrations in $\H^n$ or any hyperk\"ahler
manifold are considered in \cite{BrH}. We want to relate
our results to these works.

If $p=2$, from \ref{_explicit_form_Proposition_} we obtain
\begin{align*}
   \Psi^{2,2}_I =&
   \left.\frac{1}{2}\Re\left(\omega_K+\sqrt{-1}\omega_I\right)^2\right|_I^{2,2}
  \\ =& \left.\left(-\frac{1}{2}\omega_I^2+\frac{1}{2}\omega^2_K\right)\right|_I^{2,2}
=-\frac{1}{2}\omega^2_I+\frac{1}{4}
(\omega_J^2+\omega_K^2).
\end{align*}
  In \cite{BrH} R.Bryant and R. Harvey considered the
  forms $\Psi_{\lambda,\mu,
    \nu}=\frac{\lambda}{2}\omega_I^2+
    \frac{\mu}{2}\omega_J^2+\frac{\nu}{2}\omega_K^2$ and showed that they are calibrations iff
$-1\leq \nu,\lambda,\mu\leq1$ and $-1\leq\nu+\lambda+\mu\leq1$. We show here that the "generic" form of this type calibrates either quaternionic or complex isotropic subspaces.

\hfill

\proposition \label{_line_combi_cali_4-fo_Proposition_}
For the forms $\Psi_{\lambda,\mu,\nu}$ the following is valid:

i) If $\lambda, \mu, \nu\geq 0$ and $\lambda+\mu+\nu=1$ with at least two of $\lambda, \mu, \nu$ non-zero, the form $\Psi_{\lambda,\mu,
    \nu}$ has comass 1 and the faces are the quaternionic ones.

ii) If $\mu, \nu \leq 0$ and $\mu+\nu \geq -1$ with at least two of the inequalities being strict, then $\Psi_{1,\mu,\nu}$ has comass 1 and the faces are the I-complex
$\Omega_I$-isotropic subspaces of $\H^n=\C^{2n}$.

\hfill

{\bf Proof:} First we note that a convex hull of calibrations is a calibration. In case $i)$, for any unit 4-vector $\psi$,
\[ \Psi_{\lambda,\mu,
    \nu}(\psi)= \frac{\lambda}{2}\omega_I^2(\psi)+
    \frac{\mu}{2}\omega_J^2(\psi)+\frac{\nu}{2}\omega_K^2(\psi)
    \leq (\lambda+\mu+\nu)|\psi|=|\psi|,
\]
and the equality is achieved only when $\psi$ spans a
subspace which is invariant with respect to at least two
of $I,J$ and $K$, hence quaternionic.

    For $ii)$ we note that
\[
\frac{1}{2}\omega_I^2+
\frac{\mu}{2}\omega_J^2+\frac{\nu}{2}\omega_K^2 = \frac{1+\mu+\nu}{2}\omega_I^2-
\frac{\mu}{2}(\omega_I^2-\omega_J^2)-\frac{\nu}{2}(\omega_I^2-\omega_K^2)
\]

     Then according to \cite{BrH}, Theorem 2.38, $\frac{1}{2}(\omega_I^2-\omega_J^2)$ and $\frac{1}{2}(\omega_I^2-\omega_K^2)$ are calibrations with comass 1 and faces
which are $\omega_K$ or $\omega_J$ isotropic and contained in 2-dimensional quaternionic subspaces. So as in $i)$ if $\psi$ is a unit 4-vector, then
$\Psi_{1,\mu,\nu}(\psi)\leq|\psi|$ with equality if and only if $\psi$ is a face for all terms with nonvanishing coefficients on the right-hand-side above. If the span of $\psi$ satisfies at least two of the following:
\vspace{.005in}

i) $\psi$ is I-complex

ii) $\psi$ is $\omega_J$ isotropic and

iii) $\psi$ is $\omega_K$ isotropic
\vspace{.005in}

then $\psi$ satisfies also the third one and the Proposition follows.

\endproof

\hfill

In \cite[Theorem 6.4]{BrH},
\ref{_line_combi_cali_4-fo_Proposition_} is implicit. We
note also that in String Theory, the holomorphic
Lagrangian submanifolds in 8-dimensional manifolds were
related to the notion of intersecting branes  \cite{JF}.

\subsection{Examples}

Examples of complex Lagrangian submanifolds in
hyper-K\"ahler manifolds are given by many authors. In
\cite{_Voisin:Lagrangian_}, C. Voisin has proven a result
about the stability of such submanifolds under small
deformation of the complex structure of the ambient space;
she gave also several classes of examples.
  N. Hitchin noticed the fact that such subspaces are
 coming in complete families (\cite{Hit1}).
In \cite{_Matsushita:Lagrangian_},
D. Matsushita has shown that the families of
holomorphic Lagrangian fibrations on a hyperkaehler
manifold always deform with a deformation of a manifold,
if the cohomology class of a fiber remains of Hodge type
$(n,n)$. Existence of such families is postulated by
a conjecture called ``SYZ conjecture'', or, sometimes,
the ``Huybrechts-Sawon conjecture''. It is also known
as a hyperk\"ahler version of {\em abundance conjecture,}
related to the minimal model program.
For a survey of related questions, please see
\cite{_Sawon_}. Recently in String Theory the holomorphic
Lagrangian submanifolds were related to 3-dimensional
topological field theory with target hyperk\"ahler
manifold \cite{KRS}.

 In this section we provide examples of complex Lagrangian
 submanifolds of hypercomplex manifolds with holonomy
 $SL(n,{\Bbb H})$.

The known examples of manifolds with holonomy $SL(n,{\Bbb
  H})$ are either nilmanifolds
(\cite{_BDV:nilmanifolds_}) or
obtained via the twist construction of A.~Swann \cite{_Swann1_},
which is based on previous examples by D.~Joyce. The later
construction provides also simply-connected examples.  We describe briefly a simplified version of it.

Let $(X, I,J,K,g)$ be a compact hyper-K\"ahler
manifold. By definition, an anti-self-dual 2-form on it is
a form which is of type (1,1) with respect to $I$ and $J$
and hence with respect to all complex structures of the
hypercomplex family. Let $\alpha_1,...,\alpha_{4k}$ be
anti-self-dual closed 2-forms representing integral
cohomology classes on $X$ (instatons). Consider the principal $T^{4k}$-bundle $\pi:M\rightarrow X$ with
characteristic classes determined by $\alpha_1,...,\alpha_{4k}$. It admits an instanton connection $A$ given by
$4k$ 1-forms $\theta_i$ s.t. $d\theta_i=\pi^*(\alpha_i)$. Then $M$ carries a hypercomplex structure determined
in the following way: on the horizontal spaces of $A$ we have the pull-backs of $I,J,K$ and on the vertical
spaces we fix a linear hypercomplex structure of the
$4k$-torus. The structures $\cal{I}, \cal{J}, \cal{K}$
on $M$ are extended to act on the cotangent bundle $T^*M$
using the following relations:
\begin{align*}
{\cal I}(\theta_{4i+1})=\theta_{4i+2}, &{\cal I}(\theta_{4i+3})=\theta_{4i+4}, &{\cal J}(\theta_{4i+1})=\theta_{4i+3}, &{\cal J}(\theta_{4i+2})=-\theta_{4i+4},\\
 {\cal I}(\pi^*\alpha)=\pi^*(I\alpha), &{\cal J}(\pi^*\alpha)=\pi^*(J\alpha)& &
\end{align*}
for any 1-form $\alpha$ on $X$ and $i=0,1,...k-1$. Similarly one can define a hyperhermitian (or
quaternion-Hermitian) metric on $M$ from $g$ and a fixed hyper-K\"ahler metric on $T^{4k}$ using the splitting
of $TM$ in horizontal and vertical subspaces. As A.~Swann \cite{_Swann1_} has shown the structure is HKT and has
a holonomy $SL(n,{\Bbb H})$.

Suppose now that $Y$ is a complex Lagrangian subspace in $X$ with respect to $I$.  Consider the $T^{2k}$-bundle
over $X$ determined by $\alpha_{4i+1},\alpha_{4i+3}$. Suppose that $N$ is its restriction to $Y$ i.e $N$ is a
principal $T^{2k}$-subbundle of $M$ over $Y$ determined by $\alpha_{4i+1}|_Y,\alpha_{4i+3}|_Y$. Then $N$ is naturally
embedded in $M$ and by the definiton above $N$ is $\cal{J}$-invariant and Lagrangian with respect to the
fundamental 2-form of $\cal{I}$. Notice that in general the complex Lagrangian subspace could be K\"ahler or
non-K\"ahler depending on whether $\alpha_1|Y$ and $\alpha_3|Y$ define zero classes or not.

 As a particular case assume $X$ to be a $K3$ surface with
 large enough Picard group such that there are 4
 independent anti-self-dual integral classes defining a
 principal $T^4$-bundle $M$ over $X=K3$ with finite
 fundamental group. After passing to a finite cover we may
 assume that $M$ is simply-connected.  Now if $vol$
 denotes the volume form on $X$, then we can choose
 representatives $\alpha_1,...,\alpha_4$
in the characteristic classes of $M$ such that
 $\alpha_i^2=-F \Vol$ where $F$ is a function and $F > 0$
 almost everywhere. We want to see what is the structure
 of an arbitrary complex Lagrangian subspace $N$ of
 $M$. Since $N$ is 4-dimensional and ${\cal J}$-complex,
 we claim that its intersection with a generic fiber of
 $\pi:M\arrow X$ is at least complex
 1-dimensional. Indeed, otherwise $N$ would be a
 multisection of $M$ and will intersect a generic fiber
 transversally. However then
$\int_N \pi^*(\alpha_1^2)<0$ since $\alpha_1^2=-vol$ on
 one hand, and $\int_N \pi^*(\alpha_1^2)=0$ since
 $\pi^*(\alpha_1)=d\theta_1$ for some connection form
 $\theta_1$ on the other. The contradiction proves the
 claim and we have:

\hfill

\proposition
If $M$ is a principal instanton $T^4$-bundle over a $K3$
surface then any complex Lagrangian subspace is fibered
by complex Lagrangian curves of the fibers of $M$ over a
Lagrangian curve of the base $K3$. \endproof

\hfill

\remark
Notice that any complex curve is {\em a priori} Lagrangian in a
K3 surface.

\hfill

In general one can use a similar construction to obtain
complex isotropic and coisotropic subspaces of the
instanton bundle $M$.

%


\section{Calibrations on $SL(n, {\Bbb H})$-manifolds}


Let $(M, I, J, K, \Phi_I)$ be an $SL(n, {\Bbb
  H})$-manifold, that is, a hypercomplex manifold
with $\Phi_I$ a holomorphic volume form on $(M,I)$
preserved by the Obata connection. Clearly,
$\bar \Phi_I$ is proportional to $J(\Phi_I)$.
After a rescaling to
$e^{\sqrt{-1}t}\Phi_I$ if necessary, we can assume that
$\Phi_I$ is ${\Bbb H}$-real, i.e. $J(\Phi_I)=\bar\Phi_I$, and
${\Bbb H}$-positive (Subsection
\ref{_posi_2,0-forms_Subsection_}). A number of
interesting calibrations can be constructed in this situation.

\hfill

\theorem\label{_Lagra_calibra_on_SL_n_H_Theorem_}
Let $(M, I, J, K, \Phi_I)$ be an $SL(n, {\Bbb
  H})$-manifold, and $(\Phi_I)_J^{n,n}$ the $(n, n)$-part
of $\Phi_I$ taken with respect to $J$. Pick a quaternionic
Hermitian metric on $M$. Using a conformal change, we
may assume that $|\Phi_I|_g= 2^n$. Then $Re((\Phi_I)_J^{n,n})$ is a
calibration, and it calibrates complex subvarieties of
$(M,J)$ which are Lagrangian with respect to the
$(2,0)$-form $\omega_K + \1 \omega_I$.

\hfill

{\bf Proof:} It follows from the assumptions of
\ref{_Lagra_calibra_on_SL_n_H_Theorem_} that
\[ \Phi_I=\lambda \frac{(\omega_J+\sqrt{-1}\omega_K)^n}{n!}.\]
Since both forms are real and ${\Bbb H}$-positive,  $\lambda$ is real and
positive. It is easy to check that in local quaternionic
Hermitian frame  $(dz_1,dw_1,...,dz_n,dw_n)$ the norm is
calculated as
\[
\left|\frac{(\omega_J+\sqrt{-1}\omega_K)^n}{n!}\right|^2=
|dz_1|^2|dw_1|^2...|dz_n|^2|dw_n|^2=4^n.
\]
Then
$\left|\frac{(\omega_J+\sqrt{-1}\omega_K)^n}{n!}\right|=2^n=|\Phi_I|$
and $\lambda=1$. Now the proof follows from the fact that
$Re(\Phi_I)$ and $Re(\Phi_I)_J^{n,n}$ are
both closed,\footnote{The form $(\Phi_I)_J^{n,n}$ is
  parallel with respect to the Obata connection, which is torsion-free.}
and \ref{_complex_Lagrangian_subspaces_Proposition_}.
 \endproof

\hfill

\theorem\label{_coisotro_calibra_on_balanced_Theorem_}
Let $(M, I, J, K, \Phi_I)$ be an $SL(n, {\Bbb
  H})$-manifold, and $(\Phi_I)_J^{n,n}$ the $(n, n)$-part
of $\Phi_I$ taken with respect to $J$. Assume that $(M,I,J,K)$ is equipped with an HKT metric $g$ which is
balanced and $|\Phi_I|=2^n$. Then
$V_{n+i,n+i}:=\frac{1}{2^ii!}Re((\Phi_I)_J^{n,n}\wedge
\omega_J^i)$ is a calibration, which calibrates complex
subvarieties of $(M,J)$ which are coisotropic with respect
to the $(2,0)$-form $\omega_K + \1 \omega_I$.

\hfill

{\bf Proof:} As in the previous proof,
$\Phi_I=\frac{(\omega_J+\sqrt{-1}\omega_K)^n}{n!}$, so the form $V_{n+i,n+i}$
is a pre-calibration by \ref{_Phi^p,p_comass_Proposition_}. It is closed,
as follows from \ref{_Pi_+_of_omega^n+p_clo_Proposition_}.

 \endproof

\hfill

\remark Notice that the form $V_{n+i,n+i}$ is not
parallel with respect to any torsion-free connection
on $M$ (\ref{_conne_not_exist_Claim_}), unless $M$ is hyperk\"ahler.

\hfill

Existence of a balanced HKT metric is a hard problem, which is equivalent to a quaternionic version of a
Calabi-Yau theorem (\cite{_Verbitsky:balanced_}). However, even if $g$ is not balanced, an analogue of the
calibration $V_{n+i,n+i}$ is possible to construct.

\hfill

\theorem\label{_coisotro_calibra_on_HKT_Theorem_}
Let $(M, I, J, K, \Phi_I)$ be an $SL(n, {\Bbb
  H})$-manifold, and $(\Phi_I)_J^{n,n}$ the $(n, n)$-part
of $\Phi_I$ taken with respect to $J$, and $g$ an HKT metric. Then there exists a function $c_i(m)$ on $M$,
 such that
$V_{n+i,n+i}:=(\Phi_I)_J^{n,n}\wedge \omega_J^i$ is
a calibration with respect to the conformal metric $\widetilde{g}=c_i g$, calibrating
complex subvarieties of $(M,J)$ which are coisotropic with
respect to the $(2,0)$-form $\widetilde{\omega}_K + \1 \widetilde{\omega}_I$.

\hfill

{\bf Proof:} Since  $\Phi_I$ is ${\Bbb H}$-positive
and Obata parallel, the form $(\Phi_I)^{n,n}_J$ is closed.
Then \ref{_Pi_+_of_omega^n+p_clo_Proposition_} implies
that $V_{n+i,n+i}$ is also closed. If
we denote by $\widetilde{\omega_J}$ and
$\widetilde{\Omega_I}$ the corresponding forms after the conformal
change $\widetilde{g}=c_i(m)g$, then we can find the function $c_i(m)$ such that
\[
V_{n+i,n+i}=
\frac{1}{2^in!i!}(\widetilde{\Omega}_I^n)^{n,n}_J\wedge\widetilde{\omega}_J^i.
\]
\ref{_coisotro_calibra_on_HKT_Theorem_} then
follows from \ref{_Phi^p,p_comass_Proposition_}.
\endproof

\hfill

\remark
Similarly to the hyperk\"ahler case, it is a
natural question to ask whether the complex isotropic
submanifolds are also calibrated in $SL(n,{\Bbb
  H})$-manifolds with an HKT structure. However we can see
in the examples from Section 4.6 that this is not the
case. Consider again a toric bundle $M$ over $K3$-surface
which has 4-dimensional fiber and is
simply-connected. Such fiber contains a 2-torus which will
be a complex isotropic curve with respect to some of the
structures. By a spectral sequence argument as in Lemma
4.7 of \cite{_Swann1_}, one can see that all second
cohomology classes of $M$ are pull-backs from classes on
the base $K3$-surface. Then such a torus is homologous
to zero, since the integral of any closed 2-form on
it vanishes. Therefore, it can not be calibrated by any form.

\hfill

\claim\label{_conne_not_exist_Claim_}
Let $M$ be an $SL(n,{\Bbb H})$-manifold,
$\Omega$ an HKT-form, and $V_{n+i,n+i}$
the corresponding calibration, constructed above.
Assume that $\Omega$ is not hyperk\"ahler.
Then, the form $V_{n+i,n+i}$
is not preserved by any torsion-free connection,
for any $0<i<n$.

\hfill

{\bf Proof:} It is easy to check that
the stabilizer $St_{GL(4n, \R)}(V_{n+i,n+i})$ is equal
to the group $Sp(n)$ of quaternionic Hermitian matrices.
Therefore, any connection preserving $V_{n+i,n+i}$ would
also preserve an $Sp(n)$-structure. However, a
torsion-free connection preserving $Sp(n)$-structure
is hyperk\"ahler.
\endproof

\hfill

{\bf Acknowledgements:} We are grateful to the referee for the careful reading and many suggestions which improved the presentation of the paper.

\hfill

\hfill

\hfill

\small{

\noindent

{\sc Misha Verbitsky}\\
{\sc  Laboratory of Algebraic Geometry, SU-HSE,\\
7 Vavilova Str. Moscow, Russia, 117312}\\
{\tt verbit@maths.gla.ac.uk, \ \  verbit@mccme.ru\\}

{\sc Gueo Grantcharov\\
{\sc Department of Mathematics and Statistics\\
Florida International University\\
Miami Florida, 33199, USA}\\
\tt grantchg@fiu.edu}

 }


{\small
\begin{thebibliography}{AV1}

\bibitem[A]{_Alesker:non-commu_}
Semyon Alesker, {\em Non-commutative linear algebra and plurisubharmonic functions of quaternionic variables},
 Bull. Sci. Math. 127 (2003), no. 1, 1--35,
arXiv:math/0104209.


\bibitem[AV1]{_Alesker_Verbitsky_HKT_}
Semyon Alesker, Misha Verbitsky, {\em Plurisubharmonic functions on hypercomplex manifolds and
 HKT-geometry,} arXiv:math/0510140,
J. Geom. Anal. 16 (2006), no. 3, 375--399.


\bibitem[AV2]{_AV:Calabi_}
Semyon Alesker, Misha Verbitsky {\em Quaternionic Monge-Amp\`ere equation and Calabi problem for HKT-manifolds},
arXiv:0802.4202, Israel J. Math. 176 (2010), 109–-138.


\bibitem[BDV]{_BDV:nilmanifolds_}
Maria L. Barberis, Isabel G. Dotti, Misha Verbitsky, {\em Canonical bundles of complex nilmanifolds, with
applications to hypercomplex geometry}, arXiv:0712.3863, Math. Res. Lett. 16 (2009), no. 2, 331–-347.

\bibitem[Ber]{Berger}
M. Berger, {\em Du c\^ot\'e de chez Pu}, Ann. Scient. \'Ec. Norm. Sup. 5 (1972) no 1, 1 - 44.


\bibitem[Bes]{_Besse:Einst_Manifo_}
Besse, A., {\em Einstein Manifolds}, Springer-Verlag, New York (1987)

\bibitem[BrH]{BrH} R. Bryant, R. Harvey,
{\em Submanifolds in hyperk\"ahler geometry}, J. Am. Math. Soc., 2 (1989), no 1, 1 - 31.


\bibitem[CS]{_Capria-Salamon_}
Capria, M. M., Salamon, S. M. {\it Yang-Mills fields on quaternionic spaces}, Nonlinearity {\bf 1} (1988), no.
4, 517--530.

\bibitem[DHM]{_DHM:R^8_}
J. Dadok, R. Harvey, F. Morgan, {\em Calibrations in $\R^8$},
Trans. Amer Math. Soc., 307 (1988), 1-40.

\bibitem[F]{JF} J. Figueroa-O'Farrill
{\em Intersecting brane geometries}, J. Geom. Phys.{\bf 35} (2000), no. 2-3,



\bibitem[FG]{_Fino_Gra_}
Fino, A.,  Grantcharov, G., {\em On some properties of the manifolds with skew-symmetric torsion and holonomy
SU(n) and Sp(n)}, math.DG/0302358, Adv. Math. 189 (2004), no. 2, 439--450.



\bibitem[GP]{_Gra_Poon_}
Grantcharov, G., Poon, Y. S., {\em Geometry of hyper-K\"ahler connections with torsion},
 math.DG/9908015,
Comm. Math. Phys. 213 (2000), no. 1, 19--37.




\bibitem[HL]{_Harvey_Lawson:Calibrated_}
R. Harvey, B. Lawson, {\em Calibrated geometries},  Acta Math.  148 (1982), 47-157.


\bibitem[Hit]{Hit1}
N. Hitchin, {\em The moduli space of complex Lagrangian
  submanifolds},
Asian J. Math 3 (1999) 77-91, arXiv:math/9901069.

\bibitem[HP]{_Howe_Papado_}
P.S. Howe, G. Papadopoulos,  {\em Twistor spaces for hyper-K\"ahler manifolds with torsion} Phys. Lett. B 379
(1996), no. 1-4, 80--86.

\bibitem[H1]{_Huybrechts1_}
D. Huybrechts, {\em Compact hyperkähler manifolds}, in Calabi–Yau manifolds and related geometries (Nordfjordeid, 2001) Universitext, Springer, Berlin, 2003, 161–225.

\bibitem[H2]{_Huybrechts2_}
D. Huybrechts, {\em Hyperkähler manifolds and sheaves}, Proceedings of the International Congress of Mathematicians. Volume II, Hindustan Book Agency, New Delhi, 2010, 450–460 .

\bibitem[Hu]{_Humphreys_}
Humphreys, J., {\em Introduction to Lie Algebras and Representation
  Theory},  Graduate Texts
in Mathematics, Springer-Verlag, no. 9, 1972.


\bibitem[IP]{_Ivanov_Petkov_}
S. Ivanov, A. Petkov, {\em HKT manifolds with holonomy SL(n,H)},  arXiv:1010.5052, to appear in IMRN, 14 pages.


\bibitem[J1]{_Joyce_}
D. Joyce,  {\em Compact hypercomplex and quaternionic manifolds}, J. Differential Geom. {\bf 35} (1992) no. 3,
743-761

\bibitem[J2]{_Joyce:Calibrated_} D. Joyce, {\em Riemannian holonomy groups and calibrated geometry}, Oxford Graduate Texts in Mathematics, {\bf 12} Oxford University Press, Oxford, (2007).

\bibitem[KRS]{KRS} A.Kapustin, L. Rosansky, N. Saulina {\em Three dimensional topological field theory and symplectic algebraic geometry I}, Nucl Phys. B {\bf 816} (2009), no. 3, 295?355.

\bibitem[M]{_Matsushita:Lagrangian_}
Matsushita, D. {\em On deformations of Lagrangian fibrations},
 arXiv:0903.2098.


\bibitem[McL]{_McLean:SpLag_}
McLean, R.C. {\em Deformations of calibrated submanifolds,} Comm. Anal. Geom. 6, (1998), 705-747.

\bibitem[MS]{_Merkulov_Sch:long_}
S. Merkulov, L. Schwachh\"ofer, {\em  Classification of irreducible holonomies of torsion-free affine
connections}, math.DG/9907206, Ann. of Math. (2) 150 (1999), no. 1, 77-149, also see Addendum: math.DG/9911266,
Ann. of Math. (2) 150 (1999), no. 3, 1177-1179

\bibitem[Ob]{_Obata_}
Obata, M., {\em Affine connections on manifolds with almost complex, quaternionic or Hermitian structure}, Jap.
J. Math., 26 (1955), 43-79.






\bibitem[Saw]{_Sawon_}
Sawon, J.
{\em Abelian fibred holomorphic symplectic
  manifolds}, Turkish
Jour. Math. 27 (2003), no. 1, 197-230, math.AG/0404362.

\bibitem[Sol]{_Soldatenkov:SU(2)_}
Andrey Soldatenkov,
{\em Holonomy of the Obata connection on $SU(3)$},
arXiv:1104.2085, to appear in IMRN, 17 pages.

\bibitem[SSTV]{_Strominger:Bismut_}
A. Strominger, {\em Superstrings with torsion,} Nuclear Phys. B 274 (1986), no. 2, 253--284.

\bibitem[S]{_Swann1_}
Swann, A.,
{\em Twisting Hermitian and hypercomplex geometries},
 Duke Math. J. {\bf 155} (2010) 403 - 432, arXiv:0812.2780.

\bibitem[V0]{_so(5)_}
M. Verbitsky, {\em On the action of the Lie algebra
  $\frak{s}\frak{o}(5)$ on the cohomology of a
  hyperk\"ahler
manifold}, Func. Anal. and Appl. {\bf 24} (1990), 70-71.


\bibitem[V1]{_Verbitsky:trianalyt_}
M. Verbitsky, {\em Hyperk\"ahler embeddings and
  holomorphic symplectic geometry II}, alg-geom 9403006,
GAFA {\bf 5} no. 1 (1995), pp. 92--104.


\bibitem[V2]{_Verbitsky:desing_}
M. Verbitsky, {\em Hypercomplex Varieties}, alg-geom
9703016, Comm. Anal. Geom. 7 (1999), no. 2, 355--396.


\bibitem[V3]{_Verbitsky:HKT_}
Verbitsky, M., {\em Hyperk\"ahler manifolds with torsion,
  supersymmetry and Hodge theory}, math.AG/0112215,
Asian J. Math. Vol. 6, No. 4, pp. 679-712 (2002).


\bibitem[V4]{_Verbitsky:qD_}
M. Verbitsky, {\em Quaternionic Dolbeault complex and
  vanishing theorems on hyperkahler manifolds},
 Compos. Math. 143 (2007), no. 6, 1576--1592, math/0604303.


\bibitem[V5]{_Verbitsky:canoni_}
M. Verbitsky, {\em Hypercomplex manifolds with trivial canonical bundle
 and their holonomy}, arXiv:math/0406537,
``Moscow Seminar on Mathematical Physics, II'', American Mathematical Society Translations, {\bf 2}, 221 (2007).


\bibitem[V6]{_Verbitsky:skoda.tex_}
M. Verbitsky, {\em Positive forms on hyperkahler
  manifolds}, arXiv:0801.1899,
Osaka J. Math. Volume 47, Number 2 (2010), 353-384.


\bibitem[V7]{_Verbitsky:balanced_}
Verbitsky, M., {\em Balanced HKT metrics and strong HKT
  metrics on hypercomplex manifolds},  arXiv:0808.3218,
 Math. Res. Lett.  16  (2009),  no. 4, 735--752.





\bibitem[Vo]{_Voisin:Lagrangian_}
 C. Voisin, {\em Sur la stabilite des sous-varietes
 lagrangiennes des varietes symplectiques holomorphes}, in
 Complex projective  geometry (Trieste, 1989/Bergen,
 1989), London Math. Soc. Lecture Note Ser., 179, (1992),
 294-303.

\end{thebibliography}
}
\end{document}